\documentclass[10pt, twoside, final]{article}
\usepackage{bcc24book,url,booktabs}

\usepackage[notcite,notref,color]{showkeys} 
\usepackage{calc}
\usepackage[colorlinks,bookmarks,linkcolor=black,citecolor=black]{hyperref} 
\usepackage{enumitem}

\newenvironment{conjecture}[1][]{\ifx&#1&\begin{conj}\else\begin{conj}[#1]\fi}{\end{conj}}

\newenvironment{proposition}[1][]{\ifx&#1&\begin{prop}\else\begin{prop}[#1]\fi}{\end{prop}}
\newenvironment{question}[1][]{\ifx&#1&\begin{qn}\else\begin{qn}[#1]\fi}{\end{qn}}

\newtheorem{problem}[example]{Problem}

\setlist{itemsep=1pt,parsep=0pt,topsep=2pt,partopsep=0pt}  
\setenumerate{leftmargin=*,labelindent=\parindent} 

\def\itm#1{\rm ({#1})} 
\def\itmit#1{\itm{\it #1\,}} 
\def\rom{\itmit{\roman{*}}} 
\def\abc{\itmit{\alph{*}}}

\newcommand{\cA}{\mathcal{A}}
\newcommand{\cC}{\mathcal{C}}
\newcommand{\cG}{\mathcal{G}}
\newcommand{\cH}{\mathcal{H}}
\newcommand{\cP}{\mathcal{P}}
\newcommand{\cT}{\mathcal{T}}
\newcommand{\cX}{\mathcal{X}}
\newcommand{\cV}{\mathcal{V}}
\newcommand{\tcX}{\tilde\cX}
\newcommand{\tX}{\tilde X}

\DeclareMathOperator{\bw}{bw}
\DeclareMathOperator{\thr}{th}
\DeclareMathOperator{\Aut}{Aut}

\newcommand{\eps}{\varepsilon}
\newcommand{\Exp}{\mathbb{E}}
\newcommand{\Prob}{\mathbb{P}}
\newcommand{\dcup}{\dot\cup}
\newcommand{\mRiordan}{m_{\text{R}}}
\newcommand{\COV}{COV}
\newcommand{\LCOV}{LCOV}

\renewcommand{\subset}{\subseteq}

\bcctitle{Large-scale structures in random graphs}
\bccshorttitle{Large-scale structures in sparse random graphs}

\bccname{Julia B\"ottcher}
\bccshortname{J.~B\"ottcher}

\bccaddressa{Department of Mathematics, London School of Economics \\
  Houghton St, London WC2A 2AE, UK}
\bccemaila{j.boettcher@lse.ac.uk}

\begin{document}

\makebcctitle

\begin{abstract}
  In recent years there has been much progress in graph theory on questions
  of the following type. What is the threshold for a certain large
  substructure to appear in a random graph? When does a random graph
  contain all structures from a given family? And when does it
  contain them so robustly that even an adversary who is allowed to perturb
  the graph cannot destroy all of them? 
  I will survey this progress, and
  highlight the vital role played by some newly developed
  methods, such as the sparse regularity method, the absorbing method, and
  the container method. I will also mention many open questions that remain
  in this area.
\end{abstract}
\thispagestyle{empty}
\setcounter{page}{1}


\section{Introduction}

Erd\H{o}s and R\'enyi introduced the notion of a random graph in their
seminal paper~\cite{ErdosRenyi:first}. They thus initiated the study of
which type of property typical graphs of a certain density have or do not have,
which turned out to be immensely influential in graph theory as well as in
other related mathematical areas. The
books~\cite{Bollobas:randomGraphsBook,FriKarBook,JLRbook} provide an
excellent and extensive overview of the theory of random graphs and its
applications.

This survey is concerned with a particular type of properties of
random graphs, namely the appearance of given large-scale subgraphs. In the
past two decades, the theory of large-scale structures in random graphs $G(n,p)$
underwent swift development and originated powerful new
tools. The following three main directions of research can be distinguished in this area.

Firstly and naturally, one may study for which edge probability~$p$ the
random graph $G(n,p)$ is likely to possess \emph{one particular} spanning (or large) structure.
This structure could for example be a perfect matching, a Hamilton cycle,
or a disjoint collection of triangles covering as many vertices of $G(n,p)$
as possible.  More generally, for any sequence $(H_n)$ of graphs one can
ask when $H_n$ is a subgraph of $G(n,p)$.  Questions of this type were
pursued since the early days of the theory of random graphs, and some
turned out to be extremely challenging.

Secondly and more generally, instead of considering a single sequence
$(H_n)$ of subgraphs one may ask for $G(n,p)$ to be \emph{universal} for a
a given sequence of families $(\cH_n)$ of graphs, that is, to
simultaneously contain a copy of each graph in $\cH_n$.  A typical example
of a question of this type is for which~$p$ the random graph $G(n,p)$ is
likely to contain every binary tree on~$n$ vertices. Since the number of
such trees is huge, it is clear that an answer to this question is not
trivially entailed by a result on the appearance of any fixed spanning
binary tree in $G(n,p)$. Hence, such universality questions are in general
harder than the questions for single subgraph sequences. Universality
questions were originally motivated by problems in circuit design, data
representation, and parallel computing (see~\cite{BhaChuLeiRos} for
relevant references, and more history concerning universality). Their study
in random graphs is more recent and it is often observed or conjectured
that when $G(n,p)$ is likely to contain any fixed graph from $\cH_n$ then
it is already universal for $\cH_n$.

Finally, one may ask how \emph{resiliently} $G(n,p)$ possesses certain
structures. In other words, if $G(n,p)$ is known to contain a subgraph~$H$,
but an adversary is allowed to delete edges from $G(n,p)$ under certain
restrictions, when is the adversary likely able to destroy all copies
of~$H$. As it turns out the random graph is very robust towards such
adversarial edge deletions. Another way of motivating resilience-type
questions is from the perspective of extremal graph theory. Two main
directions of research in extremal graph theory are the investigation of
Tur\'an-type questions, and of Dirac-type questions. Tur\'an's
theorem~\cite{Turan} states that $K_r$ is a subgraph of any graph~$G$ on~$n$ vertices with more
edges than the balanced complete $(r-1)$-partite graph on~$n$ vertices
contains~$K_{r}$ as a subgraph, that is, graphs~$G$ with edge density at least
$\frac{r-2}{r-1}+o(1)$ contain $K_r$. Dirac's theorem~\cite{Dirac}, on the other hand, asserts that any
graph~$G$ with minimum degree $\delta(G)\ge\frac12 v(G)$ contains a
Hamilton cycle. Resilience-type questions then ask for the transference of
such results to sparse random graphs. For example, when does any subgraph
of $G(n,p)$ with sufficiently many edges contain~$K_r$, and when does any
subgraph of $G(n,p)$ with sufficiently high minimum degree contain a
Hamilton cycle? The former of these two questions proved to be surprisingly
deep and both questions and their generalisations inspired much recent work in the area.

In this survey I attempt to give an overview of the progress in these three main directions.
Let me stress that there is no material covered here that does not appear
elsewhere. Instead, I try to outline the exciting developments
in the area, and also give credit to the important new methods that allowed
this progress. In some cases I will give simple examples of how these
methods can be applied. These necessarily have to be brief, but pointers to
further literature will be given.

\smallskip

\noindent
{\bf What is not covered?}
There are several other important topics which recently received much
attention and are closely connected to the the developments described in
this survey in that progress in these areas influenced or was influenced by
the methods and results provided in the following, but which are, to limit
scope, not covered here. These topics include Ramsey theoretic results in
random graphs, packing results in random graphs, and embedding results in
various types of pseudorandom graphs. I also omit analogous results in
random directed graphs and random hypergraphs, and embedding results for
induced subgraphs in random graphs.

\smallskip

\noindent
{\bf Organisation.}
The survey is structured as follows.
Section~\ref{sec:def} provides basic definitions and the relevant concepts
from the theory of random graphs. Section~\ref{sec:smallH} then collects,
mainly for comparison, results on the appearance of fixed graphs~$H$ in
$G(n,p)$. Section~\ref{sec:largeH} reviews results on the appearance of a
fixed sequence $(H_n)$ in $G(n,p)$, where the graphs in $(H_n)$ grow
with~$n$, while Section~\ref{sec:universal} considers corresponding
universality results. Section~\ref{sec:robust} surveys progress on
resilience results for large subgraphs of $G(n,p)$, and
Section~\ref{sec:BUL} discusses an important tool for this type of
problem, the sparse blow-up lemma in random graphs.


\section{Basic definitions and notation}
\label{sec:def}

For easy reference, this section collects the basic definitions we need
in this survey.
Throughout, we use the natural logarithm $\log x=\log_e x$.  The
set of the first~$n$ natural numbers is denoted by $[n]=\{1,\dots,n\}$, and
$(n)_k=n\cdot(n-1)\cdot\ldots\cdot(n-k)$ is the falling factorial. As is
common in the area, ceilings and floors are omitted whenever they are not
essential.

For a graph $G=(V,E)$ we denote by~$v(G)$ the number of its vertices~$|V|$
and by~$e(G)$ the number of its edges~$|E|$. The \emph{minimum degree}
of~$G$ is $\delta(G)$, while the \emph{maximum degree} is $\Delta(G)$. The
chromatic number of~$G$ is denoted by $\chi(G)$. The \emph{girth} of a
graph is the length of its shortest cycle.
If~$H$ is a (not necessarily induced) \emph{subgraph} of~$G$ we write
$H\subset G$. An \emph{$H$-copy} in~$G$ is a (not necessarily induced) copy
of~$H$ in~$G$. The \emph{automorphism group} of~$G$ is denoted by
$\Aut(G)$.

For a vertex $v\in V$ we write $N_G(v)$ for the \emph{neighbourhood} of~$v$
in~$G$, and $\deg_G(v)=|N_G(v)|$ for its \emph{degree}. Similarly, if $U\subset V$
then $N_G(v;U)$ is the neighbourhood in~$G$ of~$v$ in the set~$U$ and
$\deg_G(v;U)=|N_G(v;U)|$. When the graph~$G$ is clear from the context we
often omit the subscript~$G$ in this notation.


\subsection{Graph classes}

The graph properties considered in this survey mainly concern the
existence of certain subgraphs, and hence are
monotone increasing. A \emph{monotone increasing graph property} is a
family~$\cP$ of graphs such that for any $G\in\cP$ we have that a
graph~$G'$ obtained from~$G$ by adding any edge is also in~$\cP$. A
monotone increasing property is \emph{non-trivial} if $K_n$ is in $\cP$ but
the complement of $K_n$ not, where $K_n$ denotes the \emph{complete graph}
on~$n$ vertices. 

A \emph{balanced $r$-partite} graph is an $r$-partite
graph whose partition classes are as equal as possible; it is \emph{complete} if
all the edges between all partition classes are present.
The cycle on~$n$ vertices is denoted by $C_n$,
and $P_n$ is  the $n$-vertex path.  A \emph{Hamilton cycle} (or \emph{path}) of
a graph~$G$ is a cycle (or path) containing all the vertices of~$G$. A
graph is called \emph{Hamiltonian} if it has a Hamilton cycle.
Let~$H$ be a fixed graph and~$G$ be a graph on~$n$ vertices. Then
an \emph{$H$-factor} in~$G$ is a collection of $\lfloor n/v(H) \rfloor$
vertex disjoint copies of~$H$. In particular, when $v(H)$ divides~$v(G)$
then an $H$-factor is a spanning subgraph of~$G$.
The \emph{$d$-dimensional cube} $Q_d$ is the graph on vertex set
$\{0,1\}^d$ with edges~$uv$ whenever~$u$ and~$v$ differ in exactly one
coordinate.  The \emph{$k\times k$-square grid} $L_k$ is the graph on
vertex set $[k]\times[k]$, with edges $uv$ whenever~$u$ and~$v$ differ in
exactly one coordinate by exactly one. 

The \emph{$k$-th power} of a
graph~$H$ is the graph obtained from~$H$ by adding all edges between
vertices of distance at most~$k$. The $2$-nd power of~$H$ is also called
the \emph{square} of~$H$. We also denote the $k$-th power of~$H$
by~$H^k$. In particular, $C_n^k$ is the $k$-th power of a cycle~$C_n$
on~$n$ vertices. A graph~$H$ is \emph{$d$-degenerate} if every subgraph
of~$H$ contains a vertex of degree at most~$d$. Equivalently, the vertices
of~$H$ can be ordered in such a way that each vertex~$v$ sends at most~$d$
edges to vertices preceding~$v$ in this order.
The \emph{bandwidth}
$\bw(H)$ of a graph~$H$ is the smallest integer~$b$ such that there is a
labelling of~$V(H)$ using all the integers~$[v(H)]$ for which $|u-v|\le b$ for each edge $uv\in
E(H)$.

Further, the following classes of graphs are considered.  Let
$\cH(n,\Delta)$ be the family of all graphs on~$n$ vertices with maximum
degree at most~$\Delta$, and $\cH(n,n,\Delta)$ be the class of all
bipartite graphs with partition classes of order~$n$ each, and with maximum
degree~$\Delta$.  The class $\cT(n,\Delta)$ contains all trees on~$n$
vertices with maximum degree~$\Delta$.
Let me remark that sometimes these graph classes will be used to refer to small
linear sized graphs, such as the class $\cH(\gamma n,\Delta)$ for some
small $\gamma>0$, where one really should write $\cH(\lfloor\gamma
n\rfloor,\Delta)$, but I omit the floors and ceilings for simplicity.


\subsection{Random graphs}

The \emph{binomial random graph} $G(n,p)$ is obtained by pairwise
independently including each of the possible $\binom{n}{2}$ edges on~$n$
vertices with probability $p=p(n)$.\footnote{This model is also often
  called the Erd\H{o}s--R\'enyi model, though this is objected to by part of
  the community because the model that Erd\H{o}s and R\'enyi used in their
  papers pioneering the area is the model $G(n,m)$.}  The \emph{uniform
  random graph} $G(n,m)$, on the other hand, assigns each graph on vertex
set $[n]$ with~$m$ edges probability $1/\binom{\binom{n}{2}}{m}$.
An event holds \emph{asymptotically almost surely} (abbreviated a.a.s.) in~$G(n,p)$
(or in $G(n,m)$) if its probability tends to~$1$ as~$n$ tends to infinity.
For a monotone increasing graph property~$\cP$ we say that $\tilde p=\tilde p(n)$ is a
\emph{threshold} for~$\cP$ if 
\[
\Prob\big(G(n,p) \in\cP\big) \to \begin{cases}
  0 \qquad & \text{if} \quad p/\tilde p \to 0 \,, \\
  1 & \text{if} \quad p/\tilde p \to \infty \,.
\end{cases}
\]
As is common, in this case $\tilde p$ will also be called \emph{the}
threshold, even though it is not unique. Bollob\'as and
Thomason~\cite{BolTho:threshold} proved that every non-trivial monotone
increasing property has a threshold. Moreover
if the threshold is of the form $\log^a n/n^b$ with $a,b>0$
fixed reals, as will be encountered frequently in this survey, then there is a
\emph{sharp} threshold~$\tilde p$, that is, for any $\eps>0$
\[
\Prob\big(G(n,p) \in\cP\big) \to \begin{cases}
  0 \qquad & \text{if} \quad p\le(1-\eps)\tilde p \,, \\
  1 & \text{if} \quad p\ge(1+\eps)\tilde p \,.
\end{cases}
\]
As explained in~\cite{Friedgut:hunting} this follows from the celebrated
characterisation of sharp thresholds by
Friedgut~\cite{Friedgut:sharp}. 

In this survey many results are considered that concern \emph{spanning}
subgraphs of~$G(n,p)$ as~$n$ tends to infinity. These results therefore do
not concern a single fixed subgraph, but rather a \emph{sequence}
of subgraphs, one for each value of~$n$. Sometimes this fact is implicitly
assumed when stating a result, but usually it is stressed by stating
that we are given a sequence $H=(H_n)$ of graphs and that $G(n,p)$
contains~$H_n$ under certain conditions a.a.s.


\subsection{Density parameters}

In the results on subgraphs~$H$ of~$G(n,p)$ that we will discuss, different
density parameters are used, which will be defined next.\footnote{It is not
  true that these parameters are densities in the sense of being
  between~$0$ and~$1$. Rather, they are variations on the average degree of
  a graph.}
The first of these parameters is called \emph{maximum
  $0$-density}, and is given by
\[m_0(H)=\max_{H'\subset H} \frac{e(H')}{v(H')}\,.\]
The maximum $0$-density is usually simply called maximum density in the
literature.  Let~$H^*$ be a subgraph of~$H$ realising the maximum in
$m_0(H)$, and let~$X^*$ be the random variable counting the number of
unlabelled copies of~$H^*$ in~$G(n,p)$. Then
\[\Exp(X^*)=\big((n)_{v(H^*)}/|\Aut(H^*)|\big)p^{e(H^*)}\approx
n^{v(H^*)}p^{e(H^*)}\,,\] which tends to infinity if $p\cdot
n^{1/m_0(H)}\to\infty$.  So, informally, we can say that in expectation
the densest subgraph of~$H$ appears around $p=n^{-1/m_0(H)}$
in~$G(n,p)$. Hence, it is natural to guess that this probability is the
threshold for the appearance of $H$-copies in $G(n,p)$, which is indeed the
case (see Theorem~\ref{thm:smallH}).

The other density parameters are slight variations on this first definition
(which can, however, have an important influence on the resulting values of
these parameters). These variations have similarly natural motivations as
the $0$-density.
The \emph{maximum $1$-density} of a graph~$H$ with at least two vertices is
\[
m_1(H)=\max_{\substack{H'\subset H \\ v(H')>1}}
\frac{e(H')}{v(H')-1}\,.
\]
This parameter is also called fractional arboricity in~\cite{AlonYuster:factor}.
Again, let~$H^*$ be a subgraph of~$H$ realising the maximum in $m_1(H)$ and
let~$v$ be a fixed vertex in~$G(n,p)$. Then the expected number of $H^*$-copies
in~$G(n,p)$ containing~$v$ tends to infinity if $p\cdot n^{1/m_1(H)}\to\infty$.
The threshold for the property that \emph{each} vertex of $G(n,p)$ is
contained in an $H$-copy is related to, but is not precisely equal to, $n^{-1/m_1(H)}$
(see also the explanations in Section~\ref{sec:smallH}).

The \emph{maximum $2$-density} of a graph~$H$ with at least one edge is
\[m_2(H)=\max_{\substack{H'\subset H \\ e(H')>0}}
\begin{cases}
  \frac{e(H')-1}{v(H')-2}, \qquad & \text{if } v(H')>2 \\
  \frac12, & \text{if } v(H')=2 \,.
\end{cases}
\]
If $p\cdot n^{1/m_2(H^*)}\to\infty$, in expectation a fixed edge of~$G(n,p)$ is contained in
many $H^*$-copies, where~$H^*$ realises the maximum in $m_2(H)$.

A graph~$H$ is called \emph{$0$-balanced} (or \emph{$1$-balanced}, or
\emph{$2$-balanced}) if~$H$ is a maximiser in $m_0(H)$ (or $m_1(H)$, or
$m_2(H)$, respectively). If~$H$ is the unique maximiser, then it is called
\emph{strictly} $0$-balanced (or $1$-balanced, or $2$-balanced, respectively).

Riordan~\cite{Riordan} defined a different density parameter $\mRiordan(H)$
for~$H$ on at least~$3$ vertices, which I will call the
\emph{maximum Riordan-density} here and which is given by
\[\mRiordan(H)=\max_{\substack{H'\subset H \\ v(H')>2}}
  \frac{e(H')}{v(H')-2}\,.
\]
Note that, again, if~$H$ is a fixed graph and the maximum in $\mRiordan(H)$
is realised by~$H^*$, then a fixed pair of vertices in~$G(n,p)$ is in
expectation contained in many $H^*$-copies if $p\cdot
n^{1/\mRiordan(H^*)}\to\infty$. Riordan, however, uses this density in a
result about copies of spanning graphs $H=(H_n)$ in~$G(n,p)$. For such a
graph~$H$, the maximum in $\mRiordan(H)$ may well be realised by some
$H^*=(H^*_n)$ with $v(H^*_n)\to\infty$, in which case $\mRiordan(H)$ is
asymptotically equal to the maximum $0$-density (or maximum $1$-density) of~$H$.


\section{Small subgraphs}
\label{sec:smallH}

This survey focuses on large subgraphs of $G(n,p)$. Before we turn to the
many results in this area, we will briefly review what is known for small,
that is, fixed subgraphs~$H$. Some of the relevant results will turn out useful
for later comparison.


\subsection{The appearance of small subgraphs}

There are a number of natural questions that one may ask concerning the
existence of a fixed subgraph~$H$ in $G(n,p)$:

\begin{enumerate}
  \item What is the threshold for the appearance of an $H$-copy?
 \item How many $H$-copies are there in $G(n,p)$?
  \item How are the $H$-copies distributed in $G(n,p)$?
\end{enumerate}

The second question is beyond the scope of this survey, though
important and strong results were obtained in this direction (see, e.g.,
\cite[Chapter~6]{JLRbook} or~\cite[Chapter~5]{FriKarBook} ). We shall concentrate on the other two,
starting with the first. A classical result by Bollob\'as~\cite{Bollobas:smallH} in the theory of
random graphs states that the threshold for the appearance of an $H$-copy
in $G(n,p)$ is determined by its maximum $0$-density.

\begin{theorem}[see, e.g., Theorem~3.4 in~\cite{JLRbook}]
\label{thm:smallH}
  Let~$H$ be a graph with (and at least one edge). The threshold for
  $G(n,p)$ to contain a copy of~$H$ is
  \[
  n^{-1/m_0(H)}\,.
  \]
\end{theorem}

This answers the first question. Note that for $0$-balanced graphs this threshold was
already established by Erd\H{o}s and R\'enyi~\cite{ErdosRenyiH}.

So let us turn to the third question, which is phrased rather vaguely. In
fact there are two meaningful interpretations which will play a more
prominent role in this survey.  One the one hand, one could ask: When do we
find many vertex disjoint copies of~$H$, or possibly even an $H$-factor?
The latter is a difficult question, and we shall return to it in
Section~\ref{sec:JKV}.  But a related question, considering a property
which is clearly necessary for an $H$-factor, is much easier: When is every
vertex of~$H$ contained in an $H$-copy? This question was answered by
Ruci\'nski~\cite{Rucinski:factor} and Spencer~\cite{Spencer:extensions}.
For strictly $1$-balanced graphs~$H$ the threshold is mainly influenced by
the maximum $1$-density of~$H$.

\begin{theorem}[see, e.g., Theorem~3.22 in~\cite{JLRbook}]
\label{thm:HatEachVertex}
Let~$H$ be a strictly $1$-balanced graph (with at least~$2$ vertices) and
let~$\COV_H$ be the event that every vertex of $G(n,p)$ is contained in a
copy of~$H$. The threshold for $\COV_H$ is
\[
\frac {(\log n)^{1/e(H)}}{n^{1/m_1(H)}} \,.
\]
\end{theorem}

Similar results for non-strictly $1$-balanced graphs exist
(see~\cite[Theorem~3.22]{JLRbook}). But these are more complicated: they
need to take into account all the different ways of rooting the graph~$H$
at some vertex and all the different subgraphs~$H'$ of~$H$ containing this
vertex. It is true, however, that the threshold for the event~$\COV_H$ of
Theorem~\ref{thm:HatEachVertex} is $\Omega(\frac{(\log
  n)^{1/e(H)}}{n^{(v(H)-1)/e(H)}})$ for every~$H$.

The appearance of the $\log$-factor in the threshold is not
surprising. Recall that the expected number of $H$-copies in $G(n,p)$
containing a fixed vertex~$v$ is of order~$n^{v(H)-1}p^{e(H)}$. Since we are
asking for an $H$ copy at \emph{every} vertex of $G(n,p)$ it is natural to
require that this quantity grows at least like $\log n$ (to allow for
concentration), which is precisely the case for $p=\frac{(\log
  n)^{1/e(H)}}{n^{1/m_1(H)}}$ if~$H$ is strictly $1$-balanced.

On the other hand, one could interpret the third question above as asking
if $G(n,p)$ has a large subgraph without any $H$-copies. It is easy to show
that this is the case below the $2$-density-threshold.

\begin{proposition}[see, e.g., Proposition~8.9 in~\cite{JLRbook}]
  For all $0<a<1$ and all~$H$ with $\Delta(H)\ge 2$ there is a constant $c>0$
  such that the following holds. If $p\le cn^{1/m_{2}(H)}$ then $G(n,p)$\ a.a.s.
  has an $H$-free subgraph~$G$ with $e(G)\ge a\cdot e(G_{n,p})$.
\end{proposition}

So $H$-copies in $G(n,p)$ are easy to delete once we are below the
$2$-density threshold. The reason for this is that the likely number of
$H$-copies is comparable to the likely number of edges at this
threshold. Above the threshold this changes, which is addressed in the
following section.


\subsection{The Erd\texorpdfstring{\H{o}}{ö}s--Stone theorem in random graphs}
\label{sec:ErdSto}

What is the maximum number of edges in an $H$-free subgraph of $G(n,p)$?
This question has inspired much research in the theory of random graphs. To
understand what the answer to this question could reasonably be, let us
first turn to dense graphs.

The Erd\H{o}s--Stone theorem, one of the cornerstones of extremal graph
theory, is a Tur\'an-type theorem which states that the crucial property of
a fixed graph~$H$ for determining the maximum number of edges in an
$H$-free graph is its chromatic number.

\begin{theorem}[Erd\H{o}s, Stone~\cite{ErdSto}]
  For each fixed graph~$H$ and every $\eps>0$ there is an~$n_0$ such that
  for all $n\ge n_0$ the following holds.  Any $n$-vertex graph~$G$ with at
  least $(\frac{\chi(H)-2}{\chi(H)-1}+\eps)\binom{n}{2}$ edges contains~$H$
  as a subgraph.
\end{theorem}

As a balanced complete $\big(\chi(H)-1\big)$-partite graph has about
$\frac{\chi(H)-2}{\chi(H)-1}\binom{n}{2}$ edges and is obviously $H$-free
this is tight up to lower order terms. Similarly, a $(\chi(H)-1)$-partite
subgraph of $G(n,p)$ with (roughly) equal sized random partition classes
and all $G(n,p)$-edges between the partition classes contains
$\big(\frac{\chi(H)-2}{\chi(H)-1}+o(1)\big)e\big(G(n,p)\big)$ edges. Hence,
when we ask for the maximum number of edges in an $H$-free subgraph of
$G(n,p)$ we cannot go below this quantity. Moreover, as explained at the
end of the last section, below the $2$-density threshold the answer becomes
(almost) trivial as all $H$-copies in $G(n,p)$ can be destroyed by deleting
just a tiny fraction of the edges.

The question then is if these two observations already tell the whole
story. The following breakthrough result on the transference of the
Erd\H{o}s--Stone theorem to sparse random graphs confirms that this is indeed
the case, and was obtained independently by Conlon and Gowers~\cite{ConGow}
(for $2$-balanced~$H$) and Schacht~\cite{Schacht} (for general~$H$).

\begin{theorem}[Schacht~\cite{Schacht}, Conlon, Gowers~\cite{ConGow}]
\label{thm:sparseErdSto}
  \mbox{} \\
  For every fixed graph~$H$ and every $\eps>0$ there are constants $0<c<C$
  such that the following holds. Let~$\cA$ be the property that the maximum
  number of edges in an $H$-free subgraph of $G(n,p)$ is at most
  $(\frac{\chi(H)-2}{\chi(H)-1}+\eps)e\big(G(n,p)\big)$. Then
  \[
  \Prob[\cA] \to \begin{cases}
    0 \qquad & \text{if } p\le cn^{-1/m_2(H)} \,,\\
    1 & \text{if } p\ge Cn^{-1/m_2(H)} \,.\\
  \end{cases}
  \]
\end{theorem}

Earlier results in this direction were obtained for special graphs~$H$ in
\cite{FraRod,Furedi,Gerke,GPSSA,GerSchSte,HaxKohLuc:evenCycles,HaxKohLuc:oddCycles,KohKreuSte,
  KohLucRod},
and for larger lower bounds on~$p$ in \cite{KohRodSch:Turan,SzaVu}.
In fact, the results of Conlon and Gowers~\cite{ConGow} and of
Schacht~\cite{Schacht} are much more general statements, allowing the
transference of a variety of extremal results on graphs, hypergraphs and sets of
integers to sparse random structures. Both proofs reduce these problems to
the analysis of random vertex subsets in certain auxiliary hypergraphs. In
the case of Theorem~\ref{thm:sparseErdSto} the vertices of the auxiliary
hypergraph~$\cH$ are the edges of~$K_n$, and the hyperedges are all
$e(H)$-tuples that form $H$-copies. Hence, the random graph $G(n,p)$
corresponds to a random subset~$S$ of $V(\cH)$, and an $H$-free subgraph of
$G(n,p)$ to an independent set in $\cH[S]$.

Recently, a very general approach has been developed to analyse such
independent sets in hypergraphs, the so-called container method developed
independently by Balogh, Morris and Samotij~\cite{BalMorSam}, and Saxton
and Thomasson~\cite{SaxTho}, which has already proved tremendously useful
for solving a variety of other problems as well.

The idea, in the language of the Erd\H{o}s--Stone theorem in random graphs,
is as follows. One naive approach to prove that $\Prob[\cA]\to 1$ if $p\ge
Cn^{-1/m_2(H)}$ in Theorem~\ref{thm:sparseErdSto} is to first fix an
$H$-free graph~$G$ on vertex set~$[n]$, to calculate the probability
that~$G(n,p)$ contains more than
$(\frac{\chi(H)-2}{\chi(H)-1}+\eps)p\binom{n}{2}$ edges of~$G$, and then
use a union bound over all choices of~$G$ to conclude that a.a.s.\ $G(n,p)$
does not contain $(\frac{\chi(H)-2}{\chi(H)-1}+\eps)p\binom{n}{2}$ edges of
any $H$-free graph, so that any subgraph of $G(n,p)$ with that many edges
cannot be $H$-free. This, of course, does not work because there are too
many choices for~$G$ (and we did not use anything about the structure of
$H$-free graphs). The crucial idea of the container method is to show that
the set~$\cG$ of $H$-free graphs can be ``approximated'' by a much smaller
set of good containers~$\cC$, that is, for each $G\in\cG$ there is
$C\in\cC$ such that $G\subset C$ and $e(C)\le
(\frac{\chi(H)-2}{\chi(H)-1}+\frac12\eps)\binom{n}{2}$. This then basically
allows us to run the union bound argument over~$\cC$. In reality, things are
not quite so simple, and more properties are required
of~$\cC$ (see also the excellent explanations
in~\cite[Section~2.2]{SaxTho}).

Restricting to the case when $H=K_r$, the analogous structural question to
Theorem~\ref{thm:sparseErdSto} of when the $K_r$-free subgraph of~$G(n,p)$
with the most edges is $(r-1)$-partite was first considered by Babai,
Simonovits and Spencer~\cite{BabSimSpe}, whose result was improved on by
Brightwell, Panagiotou and Steger~\cite{BriPanSte}. Finally, DeMarco and
Kahn~\cite{DeMarcoKahn:3,DeMarcoKahn:r} showed that this is a.a.s.\ the
case when $p\ge C(\log n)^{1/(e(K_r)-1)}/n^{m_2(K_r)}$, which is optimal up
to the value of~$C$. For other graphs~$H$ a corresponding structural result
has not yet
been established. Observe that, in general, this is a difficult problem
since we do not even know precise structural results in dense graphs for all~$H$.

\begin{question}
  For some fixed~$H$ different from a complete graph, what is the
  structure of an $H$-free subgraph of~$G(n,p)$
  with the most edges?
\end{question}

This question is already interesting when $H=C_5$, for example, when this
subgraph should be bipartite (for $p$ sufficiently large).


\section{Large subgraphs}
\label{sec:largeH}

In this section we shall consider the question when the random graph
$G(n,p)$ contains a fixed sequence of spanning graphs $H=(H_n)$ as
subgraphs. Answers to this question come in various levels of accuracy. For some
classes of graphs~$H$ we only know non-matching lower and upper bounds on
the threshold probability, while for others the threshold has been
established. Even stronger hitting time results could so far only rarely be
obtained. 

We start this section with the most classical subgraphs of~$G(n,p)$ to be
considered: matchings and Hamilton cycles. In Section~\ref{sec:Riordan} we
present a theorem of Alon and F\"uredi concerning more general spanning
subgraphs and the powerful improvement on this result by Riordan. We also
analyse the bounds on the threshold for various classes of graphs that Riordan's
result gives. Section~\ref{sec:JKV} turns to the deep Johansson--Kahn--Vu
theorem which establishes, among others, the threshold for $K_r$-factors,
while Section~\ref{sec:trees} considers the threshold for the containment
of spanning bounded degree trees. The question of when a bounded degree
graph~$H$ appears in $G(n,p)$ is addressed in Section~\ref{sec:Delta}. 
In Section~\ref{sec:KahnKalai} we discuss a very general conjecture of Kahn
and Kalai concerning the form a threshold for the containment of some sequence
$(H_n)$ can take in $G(n,p)$.
In Section~\ref{sec:constructive}, finally, we consider the question of
algorithmically finding a spanning $H$-copy in $G(n,p)$.


\subsection{Matchings and Hamilton cycles}
\label{sec:matching}

Two of the most natural questions concerning spanning substructures of random
graphs asks for the threshold of $G(n,p)$ to contain a perfect matching or
to be Hamiltonian.  These questions are as well understood as one can
hope for.

Already Erd\H{o}s and R\'enyi~\cite{ErdosRenyi:Matching_bip,
  ErdosRenyi:Matching} showed that the threshold for containing a perfect
matching is $\log n/n$.  Bollob\'as and Thomason~\cite{BolTho:smallRandom}
established a hitting time result, which considers $G(n,m)$ as a graph
process, where we start from the empty graph on~$n$ vertices and randomly
add edges one-by-one. The hitting time result then states that a.a.s.\
precisely the edge in this process which eliminates the last isolated
vertex creates a perfect matching (if~$n$ is even).  In other words,
avoiding the most trivial obstacle for containing a perfect matching in fact guarantees a perfect matching (see
also~\cite{BolFri:MatchingHamDecomp} for an alternative proof and related
results).  {\L}uczak and Ruci\'nski~\cite{LucRuc:TreeFactors} extended
these results, showing that the same hitting time result is true for
$T$-factors for any non-trivial tree~$T$.

Turning to the Hamiltonicity problem, P\'osa~\cite{Posa} and
Korshunov~\cite{Kor76, Kor77} showed that also the threshold for a Hamilton
cycle (as well as for a Hamilton path) is $\log n/n$.  Improving on this
result, Koml\'os and Szemer\'edi~\cite{KomSzem} determined an exact formula
for the probability of the existence of a Hamilton
cycle. Bollob\'as~\cite{Boll:hittingTime} established the corresponding
hitting time result, stating that as soon as $G(n,p)$ gets minimum
degree~$2$, it also contains a Hamilton cycle. Hence, if $\phi(n)$ is any
function tending to infinity, then $G(n,p)$ a.a.s.\ is Hamiltonian if
$p\ge(\log n+\log\log n+\phi(n)/n$ and not Hamiltonian if $p\le(\log
n+\log\log n-\phi(n)/n$.

Algorithmic results for finding Hamilton cycles -- a problem which is
NP-hard in general -- in random graphs above the threshold probability were
also obtained. Gurevich and Shelah~\cite{GurShe:Ham} and
Thomason~\cite{Thomason:Ham} obtained linear expected time algorithms
when~$p$ is well above the threshold.  Improving on polynomial time
randomised algorithms by Angluin and Valiant~\cite{AngVal} and
Shamir~\cite{Shamir:Ham}, Bollob\'as, Fenner and Frieze~\cite{BolFenFri}
gave a deterministic polynomial time algorithm with a success probability
that matches the probability that a Hamilton cycle exists given by Koml\'os and
Szemer\'edi~\cite{KomSzem}.


\subsection{The Alon--F\"uredi theorem and Riordan's theorem}
\label{sec:Riordan}

Let us now turn to results concerning more general
results on spanning subgraphs of $G(n,p)$.  Motivated by a question of
Bollob\'as
asking for a non-trivial probability~$p$ such that $G(n,p)$ with $n=2^d$
a.a.s.\ contains a copy of the $d$-dimensional hypercube $Q_d$, Alon and
F\"uredi~\cite{AlonFuredi} established the following result, providing an
upper bound on the threshold for the appearance of a spanning graph with a
given maximum degree. Their theorem can be seen as a first general result
concerning the appearance of spanning subgraphs in~$G(n,p)$, thus
stimulating research in the area.

\begin{theorem}[Alon, F\"uredi~\cite{AlonFuredi}]
\label{thm:AlonFuredi}
  Let $H=(H_n)$ be a fixed sequence of graphs on~$n$ vertices with
  maximum degree $\Delta(H)\le\sqrt{\sqrt{n}-1}$.
  If \[p\ge\Big(\frac{20\Delta(H)^2\log n}{n}\Big)^{1/\Delta(H)}\] then $G(n,p)$ a.a.s.\ contains a copy
  of~$H$.
\end{theorem}

Their proof uses the following simple strategy, which is based on a
multi-round exposure of $G(n,p)$.  Apply the Hajnal--Szemer\'edi
theorem~\cite{HajnalSzemeredi} to the square $H^2$ of~$H$ to obtain an
equitable $(\Delta^2+1)$-colouring of~$H^2$, that is, a partition of
$V(H^2)=V(H)$ into $\Delta^2+1$ parts $X_1\dcup\dots\dcup X_{\Delta^2+1}$
which are as equal in size as possible and form independent sets in~$H^2$.
Observe that this implies that between each pair of these parts~$H$ induces a
matching.  Partition the vertices of $G(n,p)$ into sets
$V_1\dcup\dots\dcup V_{\Delta^2+1}$ of sizes equal to these parts.  Then
embed~$X_i$ into~$V_i$ one by one, revealing the edges between $V_i$ and
$\bigcup_{j<i} V_j$, and showing that the partial embedding from the
previous round can be extended. This is possible because for any~$x\in X_i$
the set $N^-(x)$ of already embedded neighbours is of size at most~$\Delta$
and disjoint from any $N^-(x')$ with $x\neq x'\in X_i$, and a random
bipartite graph with edge probability $p^\Delta$ and partition classes of
size $n/(\Delta^2+1)$ contains a perfect matching.  Ideas from this basic
strategy were re-used in many of the results on universality and local
resilience we will mention later.

The theorem of Alon and F\"uredi was improved on by Riordan. He proved the
following surprisingly powerful result.

\begin{theorem}[Riordan's theorem~\cite{Riordan}]
\label{thm:Riordan}
  Let $H=(H_n)$ be a fixed sequence of graphs with $v(H)=n$ and $e(H)>
  n/2$ and let $p=p(n)<1$ satisfy
  \[ 
  \frac{np^{\mRiordan(H)}}{\Delta(H)^4} \to \infty \,.\]
  Then a.a.s.\ $G(n,p)$ contains a copy of~$H$.
\end{theorem}

This result can be found in this form in~\cite{ParPer} (where it is in addition
verified that this result also remains true for~$H$ with fewer edges but
$\delta(H)\ge 2$).  Observe, that the condition on~$p$ in Riordan's
theorem implies that $\Delta(H)$ grows slower than $n^{1/4}$.
In most of this survey, however, we will consider bounded degree graphs only, for
which Theorem~\ref{thm:Riordan} requires that~$p$
grows faster than $n^{-1/\mRiordan(H)}$. 

Let me mention that in~\cite{Riordan} this result is stated for~$G(n,m)$ instead
of~$G(n,p)$ and it is in addition required that $p\binom{n}{2} \to \infty$,
$\binom{n}{2}-2e(H)\to\infty$ and $(1-p)\sqrt{n}\to\infty$. However, the
result for~$G(n,p)$ follows from a standard argument (e.g.\ \cite[Theorem
2.2]{Bollobas:randomGraphsBook}) and the first additional requirement
on~$p$ follows from the requirement in Theorem~\ref{thm:Riordan}
since
$\mRiordan(H)\ge\frac{n/2}{n-2}>\frac12$ because~$e(H)>n/2$. 
The second and third additional requirements are
satisfied if we take~$p$ as small as possible while still satisfying the
conditions in Theorem~\ref{thm:Riordan} because $\Delta(H)$ grows slower
than $n^{1/4}$ and $\mRiordan(H)\le\Delta(H)$. The conclusion then still
remains true for larger~$p$ because the property of containing~$H$ is
monotone increasing.

The heart of the proof of Riordan's theorem is an elegant second moment argument
in the $G(n,m)$ model, which shows that the variance of the number of
$H$-copies is small by bounding from above how much one $H$-copy in
$G(n,m)$ can make another
$H$-copy more likely. Using the same approach in $G(n,p)$ is not possible
because if~$H$ contains many edges and one conditions on the appearance
of a fixed $H$-copy in $G(n,p)$, then this boosts the number of edges in
$G(n,p)$ sufficiently to make other $H$-copies significantly more likely.

To illustrate the power of Riordan's theorem a few straightforward
consequences are collected in the following.  The first two of these were
already given by Riordan~\cite{Riordan}, and the third was observed by
K\"uhn and Osthus~\cite{KueOst:Posa}.

\smallskip

\noindent {\bf Hypercubes.} If $n=2^d$ and \[p\ge\frac14+ 6\frac{\log
  d}{d}\] then a.a.s.\ $G(n,p)$ contains a copy of the $d$-dimensional
cube $Q_d$, because $\mRiordan(Q_d)=\frac{dn}{2(n-2)}$ 
(that is, $Q_d$ is the maximiser in $\mRiordan(Q_d)$). This
results is close to best possible since for $p=\frac14$ the expected number
of $Q_d$-copies is $(n!/|\Aut(Q_d)|)(\frac14)^{\frac12n\log
  n}\le(n!/|\Aut(Q_d)|)\cdot n^{-n}$, which tends to zero as~$n$ tends to
infinity.

\smallskip

\noindent
{\bf Square grids.} If $n=k^2$ and \[p\cdot n^{1/2}\to\infty\] then
a.a.s.\ $G(n,p)$ contains a copy of the $k\times k$-square grid~$L_k$,
because $\mRiordan(L_k)=2$ (that is, $C_4$ is the maximiser in
$\mRiordan(L_k)$).  Again, an easy first moment calculation show that for
$p=n^{-1/2}$ the probability that $G(n,p)$ contains~$L_k$ tends to~$0$.

\smallskip

\noindent {\bf Powers of Hamilton cycles.} If $k\ge 3$ and \[p\cdot
n^{1/k}\to\infty\] then~$G$ contains the $k$-th power of a Hamilton cycle
$C_n^k$, because $\mRiordan(C_n^k)\le k+\frac{(k+1)k^2}{n}$ as shown
in~\cite{KueOst:Posa}. For $p\le((1-\eps)e/n)^{1/k}$ the probability that
$G(n,p)$ contains the $k$-th power of a Hamilton cycle tends to~$0$ (using
again the first moment).

For $k=2$ Riordan's theorem does not provide a (close to) optimal result,
because $\mRiordan(C_n^k)=\mRiordan(K_3)=3$. An approximately tight result
has been obtained by K\"uhn and Osthus~\cite{KueOst:Posa} though, who
showed that $G(n,p)$ a.a.s.\ contains $C_n^2$ if $p\ge n^{\eps-1/2}$ for
any fixed $\eps>0$. This was improved on by Nenadov and
\v{S}kori\'c~\cite{NenSko:cycles} who require $p\ge C\log^4
n/n^{1/2}$. Both the result of K\"uhn and Osthus and the result of Nenadov
and \v{S}kori\'c use an absorbing-type method.  Very recently, using a
second moment argument again, Bennett, Dudek, and
Frieze~\cite{BenDudFri:square} announced a proof showing that
$p=\sqrt{1/n}$ is the threshold for $G(n,p)$ to contain the square of a
Hamilton cycle.

\smallskip

\noindent
{\bf Trees.} For trees~$T$ on at least~$3$ vertices we have
$\mRiordan(T)=2$, where the path on~$3$ vertices is the maximiser
in~$\mRiordan(T)$. It follows that if $T=(T_n)$ is a fixed sequence of
bounded degree trees then $G(n,p)$ a.a.s.\ contains~$T$ if $p \cdot
n^{1/2}\to\infty$. This is far from the best known upper bound of $\log^5
n/n$ for the threshold for containing such trees~\cite{Montgomery:DeltaT},
to which we shall return in Section~\ref{sec:trees}.
Riordan's theorem allows to also consider trees with growing maximum
degrees. However, the resulting threshold bounds are again far from the
best known bounds (see Section~\ref{sec:trees}).

\smallskip

\noindent
{\bf Planar graphs.} For a planar graph~$H'$ we have
$e(H')/(v(H')-2)\le 3$, and hence any $n$-vertex planar graph~$H$ satisfies
$\mRiordan(H)\le 3$, with equality when~$H$ is a triangulation. Hence,
if~$H$ has bounded degree, then a.a.s.\ $G(n,p)$ contains~$H$ if
\[p\cdot n^{1/3}\to\infty\,.\] As was observed by Bollob\'as and
Frieze~\cite{BolFri:planar} for $p=c/n^{1/3}$ with $c=(27e/256)^{1/3}$ the
random graph $G(n,p)$ a.a.s.\ contains no spanning triangulation.

A planar graph~$H$ drawn uniformly from all planar graphs on~$n$
vertices a.a.s.\ has maximum degree less than $3\log
n$~\cite{McDiarmidReed,DrmGimNoyPanSte}. It follows that for such
graphs~$H$ the random graph $G(n,p)$ a.a.s.\ contains~$H$ if
$p\cdot\frac{n^{1/3}}{\log^{4/3} n}\to\infty$.

\smallskip

\noindent
{\bf $K_r$-factors.} For $K_r$-factors~$H$ we have
$\mRiordan(H)=\mRiordan(K_r)=\frac12r(r-1)/(r-2)$ and hence $G(n,p)$
a.a.s.\ has a $K_r$-factor when \[p\cdot
n^{\frac2r-\frac{2}{r(r-1)}}\to\infty \,.\] The power in the exponent of~$n$
is surprisingly close to the right one, which is $-1/m_1(K_r)=-2/r$
(ignoring $\log$-factors), as given by Theorem~\ref{thm:JKV} in the next
section.

\smallskip

\noindent
{\bf Bounded degree graphs.} Graphs~$H'$ with maximum degree
$\Delta(H')\le\Delta$ satisfy
$e(H')/(v(H')-2)\le\frac12\Delta+\Delta/(v(H')-2)$. To maximise this
quantity we should set $v(H')=\Delta+1$ (since for smaller $v(H')$ an even
better bound on $e(H')/(v(H')-2)$ holds). Hence, for a maximum degree
$\Delta$ graph~$H$ we have $\mRiordan(H)\le
\frac12(\Delta+1)\Delta/(\Delta-1)$ and thus $G(n,p)$ a.a.s.\ contains~$H$
when
\[p\cdot n^{\frac{2}{\Delta+1}-\frac{2}{\Delta(\Delta+1)}}\to\infty\,.\]
This again is close to the lower bound, which is given by the lower bound
for containing a $K_{\Delta+1}$-factor.

\smallskip

\noindent
{\bf $D$-degenerate graphs.}
For $D$-degenerate graphs~$H'$ we have $e(H')\le (v(H')-D)D+\binom{D}{2}\le
v(H')D-2D$ for $D\ge 3$. It follows that a $D$-degenerate graph~$H$
satisfies $\mRiordan(H)\le D$ for $D\ge 3$. So, if further the maximum
degree of~$H$ is bounded by a constant (potentially much larger than~$D$)
then $G(n,p)$ a.a.s.\ contains~$H$ when \[p\cdot n^{1/D}\to\infty\,.\] 
As mentioned earlier for $p\le((1-\eps)e/n)^{1/D}$ the probability that
$G(n,p)$ contains the $D$-th power of a Hamilton cycle tends to~$0$. Since
the $D$-th power of a Hamilton path is $D$-degenerate this shows that the
bound given by Riordan's theorem is close to best possible.
Observe also that this bound is much better than the known
bounds in universality results for $D$-degenerate graphs discussed in
Section~\ref{sec:universalDelta}.

\smallskip

These examples illustrate that Riordan's theorem often, though not always,
gives optimal or close to optimal bounds. As indicated, for $K_r$-factors
and bounded degree trees better bounds have been obtained in recent years,
and I shall discuss these in the following sections. 

For spanning bounded degree graphs~$H$ the gap between lower bounds and the
bound given by Riordan's theorem remains, though very recently near-optimal
bounds have been obtained for almost spanning~$H$ and we shall return to
this topic in Section~\ref{sec:Delta}. 


\subsection{The Johansson--Kahn--Vu Theorem}
\label{sec:JKV}

It is not too difficult to prove (see, e.g., Theorem~4.9 of~\cite{JLRbook},
or~\cite{Rucinski:factor}) that the threshold in $G(n,p)$ for an
\emph{almost spanning}~$H$-factor, that is, a collection of vertex disjoint
copies of~$H$ covering all but at most $\eps n$ vertices, is
$n^{-1/m_1(H)}$. For obtaining a spanning $H$-factor we need to go above
this threshold by at least some (power of a) logarithmic factor
in some cases: For strictly $1$-balanced $H$, if $p$
grows slower than $(\log n)^{1/e(H)}/n^{1/m_1(H)}$ then by
Theorem~\ref{thm:HatEachVertex} a.a.s.\ not every
vertex of~$G(n,p)$ is covered by a copy of~$H$, hence $G(n,p)$ contains no
spanning $H$-factor.

Ruci\'nski~\cite{Rucinski:factor} showed that 
if $np^{\delta^*(H)}-\log n\to\infty$, where $\delta^*(H)=\max\{\delta(H')\colon
H'\subset H\}$, then $G(n,p)$ a.a.s.\ contains an
$H$-factor. This implies that the threshold for a $K_r$-factor is at
most $(\log n/n)^{1/(r-1)}$. This was improved on by
Krivelevich~\cite{Krivelevich:triangle}, who proved that for each~$r$ there
is a constant $C=C(r)$ such that if $p\ge Cn^{-2r/((r-1)(r+2))}$ then
$G(n,p)$ a.a.s.\ contains a $K_r$-factor (see~\cite[Section~4.3]{JLRbook}
for a short exposition of the interesting proof of this result in the case
$r=3$). Observe that this bound on the threshold is also better than the
one implied by Riordan's theorem (Theorem~\ref{thm:Riordan}).

Finally, in a celebrated result, Johansson, Kahn, and Vu~\cite{JKV} proved
that for strictly $1$-balanced~$H$ the threshold for an $H$-factor does
indeed coincide with the \emph{$H$-cover threshold}, that is, the threshold for
every vertex of $G(n,p)$ to be contained in an $H$-copy.

\begin{theorem}[Johansson, Kahn, Vu~\cite{JKV}] 
  \label{thm:JKV} \mbox{}\\
  For a strictly $1$-balanced graph~$H$ the threshold for $G(n,p)$ to
  contain an $H$-factor is
  \[
  \frac{(\log n)^{1/e(H)}}{n^{1/m_1(H)}} \,.
  \]
\end{theorem}

Johansson, Kahn, and Vu prove this theorem
more generally for hypergraphs in~\cite{JKV}. When~$H$ is a single
edge, that is, we are asking for a perfect hypergraph matching, it thus
solves the famous Shamir problem.  A good exposition of the proof in this
case is given in~\cite{BalFri}.

In their proof Johansson, Kahn and Vu work (for some part of the argument)
in $G(n,m)$. The basic idea is to think of $G(n,m)$ as a random graph
obtained from $K_n$ by successively \emph{deleting} random edges until
only~$m$ edges remain. They then show with the help of a martingale
argument and certain entropy results that in each deletion step not too
many $H$-factors get destroyed, implying that the number of $H$-factors in
$G(n,m)$ is close to expectation.

Already Ruci\'nski~\cite{Rucinski:factor} and Alon and
Yuster~\cite{AlonYuster:factor} observed that not for every~$H$ the
$H$-factor threshold is the same as the $H$-cover threshold. Indeed, it was
shown in~\cite{AlonYuster:factor,Rucinski:factor} that for graphs~$H$ with
$\delta(H)<m_1(H)$ the $H$-factor threshold is at least $n^{-1/m_1(H)}$, while the
$H$-cover threshold is of lower order of magnitude. It is not surprising
that the thresholds for these two properties do not always coincide since
there may be some vertex $x\in V(H)$ such that among all $H$-copies in~$G$
the vertex~$x$ is only mapped to few vertices~$u$ of~$G$.
Alon and Yuster~\cite{AlonYuster:factor} conjectured, however, that for
each graph~$H$ with $e(H)>0$ the threshold for an $H$-factor is
\[n^{-(1/m_1(H))+o(1)}\,.\]
Johansson, Kahn and Vu~\cite{JKV} prove this conjecture as well. 
Further, they conjecture that the obstacle identified in the last
paragraph is the only one, that is, that the $H$-factor threshold coincides
with the threshold for the property~$\LCOV_H$ that in an $n$-vertex
graph~$G$
\begin{enumerate}
\item each vertex of~$G$ is contained in an $H$-copy, and
\item for each $x\in V(H)$ there are at least $n/v(H)$ vertices $u\in V(G)$
  such that some $H$-copy in~$G$ maps~$x$ to~$u$.
\end{enumerate}

\begin{conjecture}[Johansson, Kahn, Vu~\cite{JKV}]
\label{conj:JKV}
  The threshold for containing an $H$-factor is the same as that for $\LCOV_H$.
\end{conjecture}

A related conjecture appears also already
in~\cite{Rucinski:factor}. Johansson, Kahn, and Vu~\cite{JKV} think it even
possible that a hitting time version of conjecture~\ref{conj:JKV} is
true. Further, they state that the threshold of $\LCOV_H$ is as follows
(for a proof see the arXiv version of~\cite[Lemma~2.5]{GerMcD}). 
The \emph{local $1$-density}
of~$H$ at $x\in V(H)$ is
\[m_1(x,H)=\max_{\substack{H'\subset H \\ x\in V(H')}} \frac{e(H')}{v(H')-1}\,.\]
We call~$H$ \emph{vertex-$1$-balanced} if $m_1(x,H)=m_1(H)$ for all $x\in
V(H)$. 
Let $s(x,H)$ denote the minimum number of edges of a maximiser~$H'$ in
$m_1(x,H)$, and let $s(H)$ be the maximum among all $s(x,H)$. 
The threshold of $\LCOV_H$, which,
following~\cite{JKV},
we denote by $\thr^{[2]}(n)$, then satisfies
\[
\thr^{[2]}(n)=\begin{cases}
  \frac{(\log n)^{1/s(H)}}{n^{1/m_1(H)}} \quad & \text{if $H$ is
    vertex-$1$-balanced} \,,\\
  {n^{-1/m_1(H)}} & \text{otherwise}\,.
\end{cases}
\]

Gerke and McDowell~\cite{GerMcD} proved Conjecture~\ref{conj:JKV} for
graphs which are not
vertex-$1$-balanced. Hence, the only open case now is that of
vertex-$1$-balanced graphs which are not strictly $1$-balanced.

\begin{theorem}[Gerke, McDowell~\cite{GerMcD}]
\label{thm:GerMcD}
  For a graph~$H$ which is not vertex-$1$-balanced the threshold for an
  $H$-factor in $G(n,p)$ is $n^{-1/m_1(H)}$.
\end{theorem}

The idea of~\cite{GerMcD} is to identify dense subgraphs $H'$ of~$H$ (which
do not cover all vertices because~$H$ is non-vertex-$1$-balanced) and first
embed a corresponding non-spanning $H'$-factor into $G(n,p)$.  They then
use a variant of Theorem~\ref{thm:JKV} to complete the embedding.  For
obtaining this variant they verify that a partite version of the
Johansson--Kahn--Vu theorem holds, which is also useful in other
applications.

In fact, the method of Gerke and McDowell allows a proof of
Conjecture~\ref{conj:JKV} also in the case of many~$H$ which are
vertex-$1$-balanced and not strictly $1$-balanced. Moreover, for all
other~$H$ (as for example a triangle and a $C_4$ glued along one edge) the
upper bound given by their method is within a constant log-power of the
conjectured bound (see the discussions in the concluding remarks
of~\cite{GerMcD}).


\subsection{Trees}
\label{sec:trees}

The appearance of long paths in~$G(n,p)$ was another topic considered early
on in the theory of random graphs.  As explained in
Section~\ref{sec:matching} the threshold in $G(n,p)$ for a Hamilton path is
$\log n/n$, where the lower bound follows from the fact that for $p<\log
n/n$ there are a.a.s.\ isolated vertices in $G(n,p)$.  Many related results
were obtained in the sequel. To give an example, in
\cite{AjtKomSze:paths,Fernandez:paths} paths of length~$cn$ in $G(n,p)$ for
$0<c<1$ are considered. But one very natural question, which turned out to
be difficult, is if the threshold result for Hamilton paths extends to
other spanning trees with bounded maximum degree.  The following
conjecture, which claims that this is indeed the case and has prompted much
recent work, is attributed to Kahn (see~\cite{KahLubWor:comb}), but also
appears in~\cite{AloKriSud:trees}.

\begin{conjecture}
\label{conj:trees}
  For every fixed~$\Delta$ there is some constant~$C$ such that if
  $T=(T_n)$ is a fixed sequence of trees on~$n$ vertices with
  $\Delta(T)\le\Delta$ then $G(n,p)$ a.a.s.\ contains~$T$ if $p\ge C\log n/n$.
\end{conjecture}

In the following I will summarise the progress that has been made towards
proving this conjecture. Trees of small linear size were considered by
Fernandez de la Vega~\cite{Fernandez:trees}, who proved that there are
(large) constants $C,C'$ such that for any fixed $\Delta$ and any fixed
sequence $T=(T_n)$ of trees with $v(T)\le n/C$ and $\Delta(T)\le\Delta$, if
$p\ge C'\Delta/n$ then $G(n,p)$ a.a.s.\ contains~$T$. Alon, Krivelevich and
Sudakov~\cite{AloKriSud:trees} improved on this and showed that the
threshold in $G(n,p)$ for any sequence of almost spanning trees of bounded
degree is $1/n$.

\begin{theorem}[Alon, Krivelevich, Sudakov~\cite{AloKriSud:trees}]
\label{thm:AKStrees}
  Given $\Delta\ge 2$ and $0<\eps<\frac12$, let
  $C=10^6\Delta^3\eps^{-1}\log\Delta\log^2(2/\eps)$. If $p\ge C/n$ then
  $G(n,p)$ a.a.s.\ contains all trees~$T$ with $\Delta(T)\le\Delta$ and $v(T)\le(1-\eps)n$.
\end{theorem}

Observe that Theorem~\ref{thm:AKStrees} is a universality result, stating
that $G(n,p)$ contains all these trees \emph{simultaneously}. We shall
discuss universality results in $G(n,p)$ in more detail in
Section~\ref{sec:universal}. Obtaining such a universality result is
possible for Alon, Krivelevich, and Sudakov because they do not prove their
result directly for $G(n,p)$, but instead for any graph satisfying certain
degree and expansion properties. Their proof uses the well-known embedding
result for small (linear sized) trees by Friedman and
Pippenger~\cite{FriedPipp}. Balogh, Csaba, Pei and
Samotij~\cite{BalCsaPeiSam:trees} showed that using instead a related
tree embedding result of Haxell~\cite{Haxell:FriedPipp}, which works for
larger trees, one can improve the constant in Theorem~\ref{thm:AKStrees} to
$C=\max\{1000\Delta\log(20
\Delta),30\Delta\eps^{-1}\log(4e\eps^{-1})\}$. This was further improved by
Montgomery~\cite{Montgomery:comb} to
$C=30\Delta\eps^{-1}\log(4e\eps^{-1}))$, which comes close to the
$C=\Theta(\Delta\log\eps^{-1})$ believed possible
in~\cite{AloKriSud:trees}.

Alon, Krivelevich and Sudakov also observed in~\cite{AloKriSud:trees} that
for every $\eps>0$ Theorem~\ref{thm:AKStrees} immediately implies
Conjecture~\ref{conj:trees} for trees~$T$ with $\eps n$ leaves (for $p\ge
C(\eps,\Delta)\log n/n$), by using a two-round exposure of $G(n,p)$,
finding in the first round a copy of $T$ minus $(\eps n/\Delta)$ leaves
with distinct parents, and then embedding these leaves in the second round,
which is easy because all it requires is to find a certain matching.
Hefetz, Krivelevich, and Szab\'o observe in~\cite{HefKriSza:trees} that a
similar strategy can be used for embedding trees~$T$ with a linearly sized
\emph{bare path}, that is a path whose inner vertices have degree~$2$
in~$T$, also for $p\ge C(\eps,\Delta)\log n/n$.

This leaves the case of trees with few leaves (and no long bare path) of
Conjecture~\ref{conj:trees}. Since each tree has average degree less
than~$2$, however, these trees have many vertices of degree~$2$, and
hence a linear number of (arbitrarily long) constant length bare
  paths. Krivelevich~\cite{Krivelevich:trees} used this fact and showed that
the same strategy as outlined for trees with many leaves in the previous
paragraph can be used for trees with many bare paths by replacing the
matching argument by a partite version of the Johansson--Kahn--Vu theorem
for embedding the bare paths. Krivelevich's strategy leads to the following result.

\begin{theorem}[Krivelevich~\cite{Krivelevich:trees}]
  \label{thm:KrivelevichTrees}
  For every $\eps>0$ and every sequence $T=(T_n)$ of trees with $v(T)\le n$
  the random graph $G(n,p)$ a.a.s.\ contains~$T$ if
  \[p\ge\frac{40\Delta(T)\eps^{-1}\log n +n^{\eps}}{n}\,.\]
\end{theorem}

In this result $\Delta(T)$ is allowed to grow
with~$n$ (in particular, a different strategy than
Theorem~\ref{thm:AKStrees} is used for obtaining an almost spanning
embedding).

Further progress on various classes of trees has been obtained by various
groups.  Hefetz, Krivelevich, and Szab\'o~\cite{HefKriSza:trees} show that
trees with linearly many leaves and trees with linear sized bare paths, and
Montgomery~\cite{Montgomery:comb} that trees with $\alpha n/\log^9 n$ bare
paths of length $\log^9 n$ for any $\alpha>0$, are already a.a.s.\
contained in $G(n,p)$ for \[p=(1+\eps)\log n/n\,.\] Hefetz, Krivelevich, and
Szab\'o~\cite{HefKriSza:trees} also argue that for
the same~$p$ the random graph $G(n,p)$
a.a.s.\ contains any typical random tree~$T$, that
is, a tree with maximum degree $(1+o(1))\log n/\log\log n$
 as shown in~\cite{Moon:DeltaT}.

Investigating a class of special trees called combs was suggested by Kahn
(see~\cite{KahLubWor:comb}). A \emph{comb} is a tree consisting of a path
on $n/k$ vertices with disjoint $k$-paths beginning at each of its
vertices. Observe that, for example for $k=\sqrt{n}$, combs neither have
linearly many leaves nor linear sized bare paths. Kahn, Lubetzky, and
Wormald~\cite{KahLubWor:comb,KahLubWor:comb2} established
Conjecture~\ref{conj:trees} for combs. This was improved on and generalised
by Montgomery~\cite{Montgomery:comb} who proved the following result.  A
\emph{tooth} of length~$k$ in a tree is a bare path of length~$k$ where one
end-vertex is a leaf. Montgomery showed that for any fixed $\alpha>0$ a
tree $T$ with at least $\alpha n/k$ teeth of length~$k$ is contained
a.a.s.\ in $G(n,p)$ for $p=(1+\eps)\log n/n$.

Finally, a result for general bounded degree trees has recently been
established by Montgomery~\cite{Montgomery:DeltaT}, which comes very close
to the conjectured threshold. 

\begin{theorem}[Montgomery~\cite{Montgomery:DeltaT}]
\label{thm:Montgomery}
  If~$T=(T_n)$ is a fixed sequence of trees on~$n$ vertices with maximum degree
  $\Delta=\Delta(n)$ then $G(n,p)$ a.a.s.\ contains~$T$ if
  $p\ge\Delta\log^5 n/n$.
\end{theorem}

Montgomery also announced in~\cite{Montgomery:DeltaT} further work in
progress leading to the proof of Conjecture~\ref{conj:trees}.  For proving
Theorem~\ref{thm:Montgomery} Montgomery follows the basic strategy outlined above of first
finding an almost spanning subtree of~$T$, leaving some bare paths to be
embedded in a second stage (since the case of trees with many leaves is
solved already). For embedding these bare paths, however, Montgomery uses
an absorbing-type method. 


\subsection{Bounded degree graphs}
\label{sec:Delta}

Now we turn to the question of when $G(n,p)$ contains given spanning graphs
of bounded maximum degree.  Let $\Delta$ be a constant and $H=(H_n)$ be
sequence of graphs with $\Delta(H)\le\Delta$ and $v(H)\le n$.  Recall that
the Theorem of Alon and F\"uredi (Theorem~\ref{thm:AlonFuredi}) implies
that $G(n,p)$ a.a.s.\ contains~$H$ if
$p\ge (C(\Delta)\log n/n)^{1/\Delta}$, and Riordan's theorem
(Theorem~\ref{thm:Riordan}) implies the same if
$p\cdot n^{\frac{2}{\Delta+1}-\frac{2}{\Delta(\Delta+1)}}\to\infty$.  This
is unlikely to be optimal, though it cannot be far off. The optimum is widely
believed to be as follows (see, e.g., \cite{FerLuhNgu:Delta}).
\begin{conjecture}
  \label{conj:Delta}
  Let $H=(H_n)$ be a sequence of graphs with $\Delta(H)\le\Delta$ and
  $v(H)\le n$. Then
  $G(n,p)$ a.a.s.\ contains~$H$ if
  \begin{equation}
    \label{eq:Delta}
    p \cdot\frac{n^{2/(\Delta+1)}}{(\log
      n)^{1/\binom{\Delta+1}{2}}}\to\infty\,.
  \end{equation}
\end{conjecture}
In other words, the conjecture states that $G(n,p)$ contains~$H$ from above
the threshold for a $K_{\Delta+1}$-factor. Ferber, Luh and
Nguyen~\cite{FerLuhNgu:Delta} prove Conjecture~\ref{conj:Delta} for almost
spanning~$H$.

\begin{theorem}[Ferber, Luh, Nguyen~\cite{FerLuhNgu:Delta}]
  \label{thm:DeltaAlmost}
  Let $\eps>0$ and~$\Delta$ be fixed. Let $H=(H_n)$ be a fixed sequence of
  graphs with $\Delta(H)\le\Delta$ and $v(H)\le(1-\eps)n$. Then $G(n,p)$
  a.a.s.\ contains~$H$ if~$p$ satisfies~\eqref{eq:Delta}.
\end{theorem}

The strategy for the proof of Theorem~\ref{thm:DeltaAlmost} is as
follows. Ferber, Luh and Nguyen show that~$H$ can be partitioned into a
sparse part~$H'$, which is sparse enough to be embedded with the help of
Riordan's theorem, and a dense part which consists of a collection of
induced subgraphs, each of constant size. Given an embedding of~$H'$ they
then in constantly many rounds extend this
embedding successively to embed also the constant size dense bits of~$H$
by finding a matching in a suitable auxiliary hypergraph, using a
hypergraph Hall-type theorem of Aharoni and Haxell~\cite{AharoniHaxell} (a
similar idea was already used in~\cite{ConFerNenSko:universal}).
                                                            
In fact, it is widely believed that even a
universality version of Conjecture~\ref{conj:Delta} is true (see
Conjecture~\ref{conj:universalDelta}). Further recent advances were made in
this direction, which we shall return to in
Section~\ref{sec:universalDelta}.


\subsection{The Kahn--Kalai conjecture}
\label{sec:KahnKalai}

Let us round off the results presented in the previous sections with a
far-reaching and appealing conjecture of Kahn and Kalai. We first need some
motivation and definitions.  Theorem~\ref{thm:smallH} states that for fixed
graphs~$H$ the threshold for the appearance of~$H$ in $G(n,p)$ coincides
with what Kahn and Kalai~\cite{KahnKalai} call the \emph{expectation threshold} for~$H$, written $p_\Exp(H,n)$, which is the least $p=p(n)$ such that
for each subgraph~$H'$ of~$H$ the expected number of~$H'$ in $G(n,p)$ is at
least~$1$. The expectation threshold can be defined analogously for
sequences $H=(H_n)$ of graphs. In particular, for any $(H_n)=H$ we have
that $p_\Exp(H,n)$ is the least $p=p(n)$ such that for every subgraph~$H'$
of~$H$ we have \[\frac{(n)_{v(H')}}{|\Aut(H')|}p^{e(H')}\ge 1 \,.\] For
example, if~$H$ is an $F$-factor and $F'$ is any subgraph of~$F$ then
let~$H'$ be the vertex disjoint union of $\ell=n/v(F)$ copies of~$F'$. Then
the condition above requires that
\[\frac{(n)_{\ell v(F')}}{\ell!|\Aut(F')|}p^{\ell e(F')}\ge
1\,,\]
which can easily be calculated to be equivalent to $p\ge
Cn^{-(v(F')-1)/e(F')}$ for some constant~$C$, and hence $p_\Exp(H,n)$ is of
the order $n^{-m_1(F)}$ for $F$-factors. So, by
Theorem~\ref{thm:JKV}, in this case $p_\Exp(H,n)$ is different from the
threshold for the appearance of~$H$ if $F$ is strictly balanced -- but only
by less than a $\log n$ factor. Kahn and Kalai~\cite{KahnKalai} conjectured
that this is the case for every~$H$.

\begin{conjecture}[Kahn, Kalai~\cite{KahnKalai}]
\label{conj:KahnKalai}
  There is a universal constant~$C$ such that for any sequence $H=(H_n)$ of
  graphs the threshold for $G(n,p)$ to contain~$H$ is at most
  $Cp_\Exp(H,n)\log n$.
\end{conjecture}

Conjecture~\ref{conj:trees} on trees is a special case of
Conjecture~\ref{conj:KahnKalai} because $p_\Exp(T,n)$ is of order $1/n$ for
bounded degree trees~$T$. Conjecture~\ref{conj:JKV} on $H$-factors and
Conjecture~\ref{conj:Delta} on bounded degree graphs, on the other hand,
are somewhat stronger than what is implied by
Conjecture~\ref{conj:KahnKalai} because they specify a smaller
$\log$-power.


\subsection{Constructive proofs}
\label{sec:constructive}

One question we have only occasionally taken up in the preceding sections
is if the results on the various structures that exist in~$G(n,p)$ a.a.s.\
for certain probabilities have constructive proofs,
allowing for a deterministic or randomised algorithm
which finds the desired structure. This question is important for two reasons:

\begin{enumerate}
  \item Such constructive proofs often lead to polynomial time algorithms,
    making it possible to find the structures efficiently.
  \item Constructive proofs often allow the identification of certain
    pseudorandom properties, that is, properties which $G(n,p)$ a.a.s.\
    enjoys, which are sufficient for the construction to
    work. In this case universality results may become possible. 
\end{enumerate}

In particular, two prominent results we discussed, whose proofs were not
constructive but used the second moment method, were Riordan's theorem and
the Johansson--Kahn--Vu theorem. As outlined, these were also used as tools
in the proof of other results, such as Theorem~\ref{thm:GerMcD},
Theorem~\ref{thm:KrivelevichTrees}, or Theorem~\ref{thm:DeltaAlmost}.
This motivates the following problem.

\begin{problem}
  Give a constructive proof of Riordan's theorem
  (Theorem~\ref{thm:Riordan}) or the Johansson--Kahn--Vu theorem
  (Theorem~\ref{thm:JKV}).
\end{problem}

As I shall explain in Sections~\ref{sec:universal} and~\ref{sec:robust},
many constructive proofs for embedding classes of spanning or almost
spanning graphs~$H$ in~$G(n,p)$ (or in subgraphs of~$G(n,p)$) we know of
follow a greedy-type paradigm: They embed~$H$ (or a suitable subgraph
of~$H$) vertex by vertex (or class of vertices by class of vertices),
aiming at guaranteeing that unembedded common $H$-neighbours of already
embedded $H$-vertices can still be embedded in the future. In this sense
they crucially rely on the fact that all common neighbourhoods in $G(n,p)$
of $\Delta(H)$ vertices (or of~$D$ vertices if~$H$ is $D$-degenerate) are
large, which fails to be true for $p\le n^{-1/\Delta(H)}$. Hence, in
Sections~\ref{sec:universal} and~\ref{sec:robust} probability bounds of
this order shall often form a natural barrier not yet overcome in many
instances, though they are not believed to be the right bounds.


\section{Universality of random graphs}
\label{sec:universal}

In this section we consider the question of when the random graph is
a.a.s.\ universal for certain classes of graphs. More precisely, a
graph~$G$ on~$n$ vertices is said to be \emph{universal} for a class~$\cH$ of
graphs, if it contains a copy of every graph $H\in\cH$.
The crucial difference for $G(n,p)$ to contain some
$H\in\cH$ a.a.s.\ and to be a.a.s.\ universal for~$\cH$ (if~$\cH$ is large)
is that in the latter case we require a typical graph from $G(n,p)$ to
contain all these $H\in\cH$ \emph{simultaneously}.

The graph classes for which universality results have been established, and
which we shall consider in this section are bounded degree graphs, bounded
degree graphs which further have (smaller) bounded maximum $0$-density or
bounded degeneracy, and bounded degree trees. Let me stress that none of
the results presented in this section is believed to be optimal,
indicating that the methods we have at hand for proving universality are
still limited. Moreover, there are many other natural graph classes still
to be considered. The following is just one example.

\begin{question}
  When is $G(n,p)$ a.a.s.\ universal for the class of all
  planar graphs with maximum degree~$\Delta$; or more generally for all
  maximum degree~$\Delta$ graphs which are $F$-minor free for some
  fixed~$F$?
\end{question}

As an aside, $n$-vertex universal graphs with $O(n\log n)$ edges for
$n$-vertex planar graphs with maximum degree~$\Delta$ were constructed
in~\cite{BhaChuLeiRos}, and graphs~$G$ with $v(G)+e(G)=O(n)$ that are
universal for this class of graphs in~\cite{Capalbo:planar}.  For more
background on constructions of universal graphs see the survey of
Alon~\cite{Alon:ChaosOrder}.

\subsection{Universality for bounded degree graphs} 
\label{sec:universalDelta}

Before we turn to results concerning the universality of $G(n,p)$ for the
family $\cH(n,\Delta)$ of all $n$-vertex graphs with maximum degree at
most~$\Delta$, let us first briefly recall some lower bounds. A counting
argument shows that any graph~$G$ that is universal for $\cH(n,\Delta)$ must
have edge density at least $\Omega(n^{-2/\Delta})$. This was observed
in~\cite{millenium}, and follows from the fact that $\sum_{i\le \Delta n/2}
\binom{e(G)}{i}\ge|\cH(\Delta,n)|$ and well-known estimates of the number
of $\Delta$-regular graphs (for details see~\cite{millenium}). It is
interesting to observe that this lower bound was matched by constructive
results: Alon and Capalbo constructed graphs that are universal for
$\cH(n,\Delta)$ and have~$n$ vertices and
$C(\Delta)n^{2-2/\Delta}\log^{4/\Delta} n$ edges
in~\cite{AlonCapalbo:universalSpanning}, and $(1+\eps)n$ vertices and
$C_2(\Delta,\eps)n^{2-2/\Delta}$ edges for every $\eps>0$
in~\cite{AlonCapalbo:universalOptimal} (see also~\cite{Alon:ChaosOrder}).

For $G(n,p)$ the only better lower bound we know is the following, which
is only slightly better and only appeals to one particular graph in~$\cH(n,\Delta)$
instead of universality. By Theorem~\ref{thm:HatEachVertex}, If $p$ grows slower than $(\log
n)^{1/\binom{\Delta+1}{2}}/n^{2/(\Delta+1)}$ then a.a.s.\ $G(n,p)$ contains
no spanning $K_{\Delta+1}$-factor. If one turns to universality for smaller, but linearly sized
graphs~$H$, the known lower bound is not much smaller. Indeed, if $p\le c
n^{-2/(\Delta+1)}$ for some sufficiently small~$c=c(\eta)>0$ then $G(n,p)$
is not universal even for $\cH(\eta n, \Delta)$ as it does not contain a
vertex disjoint union of $K_{\Delta+1}$ covering $\eta n$ vertices
because the expected number of $K_{\Delta+1}$ in $G(n,p)$ is at most
$
 n^{\Delta+1}p^{(\Delta+1)\Delta/2}
 \le c^{(\Delta+1)\Delta/2}n
$.

As mentioned earlier, it is widely believed (see, e.g., \cite{DelKohRodRuc,FerLuhNgu:Delta}) that the lower bound above reflects the truth, that is,
when $G(n,p)$ starts containing every fixed sequence $(H_n)$
of graphs from $\cH(n,\Delta)$ a.a.s.\ then it is already universal for
$\cH(n,\Delta)$ (cf.\ Conjecture~\ref{conj:Delta}). 

\begin{conjecture}
\label{conj:universalDelta}
  $G(n,p)$ is a.a.s.\ universal for $\cH(n,\Delta)$ if
  \begin{equation*}
   p \cdot\frac{n^{2/(\Delta+1)}}{(\log
      n)^{1/\binom{\Delta+1}{2}}}\to\infty\,.
  \end{equation*}
\end{conjecture}

At present we are still far from verifying
Conjecture~\ref{conj:universalDelta}, though this problem attracted
considerable attention since the turn of the millennium.  Alon, Capalbo,
Kohayakawa, R\"odl, Ruci\'nski and Szemer\'edi~\cite{millenium} considered
almost spanning graphs and showed that for every $\eps>0$ and $\Delta$
there is~$C$ such that for $p\ge C(\log n/n)^{1/\Delta}$ the random graph
$G(n,p)$ is a.a.s.\ universal for $\cH\big((1-\eps)n,\Delta\big)$.  After
improvements in~\cite{DelKohRodRuc:suboptimal}, Dellamonica, Kohayakawa,
R\"odl, and Ruci\'nski~\cite{DelKohRodRuc} showed that for this probability
$G(n,p)$ is also universal for spanning bounded degree graphs.

\begin{theorem}[Dellamonica, Kohayakawa, R\"odl, Ruci\'nski~\cite{DelKohRodRuc}]
\label{thm:DKRRuniversal}
\mbox{}\\
  For each $\Delta\ge 3$ there is~$C$ such that
  $G(n,p)$ is a.a.s.\ universal for the family $\cH(n,\Delta)$
  if
  \[p\ge C\Big(\frac{\log n}{n}\Big)^{1/\Delta}\,.\]
\end{theorem}

Using a simpler argument (but the same basic strategy), Kim and
Lee~\cite{KimLee:universal} showed that this result also holds for
$\Delta=2$.  For proving their theorem Dellamonica, Kohayakawa, R\"odl, and
Ruci\'nski present a randomised algorithm that uses a certain set of
pseudorandom properties which $G(n,p)$ has a.a.s.\ and embeds every
$H\in\cH(n,\Delta)$ a.a.s.\ in every graph~$G$ with these pseudorandom
properties. This algorithm is inspired by the various known techniques for
proving the blow-up lemma in dense
graphs~\cite{KSS:blowup,KSS:blowup_alg,RodRuc:blowup,RodRucTar:blowup},
and the underlying idea of using an embedding strategy based on matchings
goes back to the proof of the theorem of Alon and F\"uredi
(Theorem~\ref{thm:AlonFuredi}) outlined in Section~\ref{sec:Riordan}.

As mentioned in Section~\ref{sec:constructive} the exponent $1/\Delta$
forms a natural barrier to further improvement. So far, this barrier was
broken only for almost spanning subgraphs and in the case $\Delta=2$.

\begin{theorem}[Conlon, Ferber, Nenadov, \v{S}kori\'c~\cite{ConFerNenSko:universal}]
\label{thm:CFNSuniversal}
\mbox{}\\
  For every $\eps>0$ and $\Delta\ge 3$ the random graph $G(n,p)$ is a.a.s.\
  universal for $\cH\big((1-\eps)n,\Delta\big)$ if
  \[p\cdot\frac{n^{1/(\Delta-1)}}{\log^5 n}\to\infty\,.\]
\end{theorem}

For the case of maximum degree $\Delta=2$, Conlon, Ferber, Nenadov, and
\v{S}kori\'c~\cite{ConFerNenSko:universal} also state that similar
arguments as those used for showing this theorem show that $G(n,p)$ is
a.a.s.\ universal for $\cH\big((1-\eps)n,2\big)$ if $p\ge Cn^{-2/3}$, which
is best possible up to the value of~$C$. Moreover, Ferber, Kronenberg, and Luh~\cite{FerKroLuh} very
recently showed that $G(n,p)$ is a.a.s.\ universal for $\cH(n,2)$ if $p\ge
C(\log n/n^2)^{1/3}$, which is again best possible up to the value
of~$C$. Their proof combines the Johansson--Kahn--Vu Theorem with arguments
from Montgomery's~\cite{Montgomery:DeltaT} proof of
Theorem~\ref{thm:Montgomery}.

The strategy of the proof of Theorem~\ref{thm:CFNSuniversal} is as
follows. Each graph~$H$ under consideration is partitioned into a
set of (small) components with at most $\log^4 n$ vertices, a set of
induced cycles of length at most $2\log n$, and the graph~$H'$ induced on
the remaining vertices. They then show that any induced subgraph of
$G(n,p)$ on $\frac12\eps n$ vertices is universal for $\cH(\log^4
n,\Delta)$ and can thus be used for embedding the small components,
that~$H'$ has a structure suitable for a technical embedding result of
Ferber, Nenadov and Peter~\cite{FerNenPet:universal}, and that the
remaining short cycles can be embedded with the help of the hypergraph
matching criterion of Aharoni and Haxell~\cite{AharoniHaxell}.

In~\cite{FerNenPet:universal} Ferber, Nenadov and Peter use the technical
result just mentioned for a spanning universality result under additional
constraints. More precisely, they consider graphs in $\cH(n,\Delta)$ with maximum
$0$-density at most~$m_0$, and provide a better bound than
Theorem~\ref{thm:DKRRuniversal} for $m_0<\Delta/4$. Note that for any~$H$
we have $m_0(H)\le\Delta(H)/2$.

\newlength{\myspace}
\setlength{\myspace}{\widthof{$\text{ and}$}}

\begin{theorem}[Ferber, Nenadov, Peter~\cite{FerNenPet:universal}]
\label{thm:FNPuniversal}
  For $\Delta=\Delta(n)>1$ and $m_0=m_0(n)\ge 1$ the random graph $G(n,p)$
  is a.a.s.\ universal 
  \begin{enumerate}[label=\abc,labelindent=0pt]
    \item for all $H\in\cH(n,\Delta)$ with $m_0(H)\le m_0$ if
      \[p\cdot\frac{\Delta^{12} n^{1/(4m_0)}}{\log^3 n}\to\infty\,,\text{ and}\]
    \item for all $H\in\cH(n,\Delta)$ with $m_0(H)\le m_0$ and girth at
      least~$7$ if
      \[p\cdot \frac{\Delta^{12} n^{1/(2m_0)}}{\log^3 n}\to\infty\,.
      \hspace{\myspace}\]
  \end{enumerate}
\end{theorem}

Ferber, Nenadov and Peter prove this result by using a similar embedding
strategy (and a similar decomposition of the graphs~$H$) as Dellamonica,
Kohayakawa, R\"odl, Ruci\'nski~\cite{DelKohRodRuc} and Kim and
Lee~\cite{KimLee:universal}. 

A related result is proven in~\cite{sparseBlowUp}, where $D$-degenerate
graphs~$H$ in $\cH(n,\Delta)$ are considered.  It is not difficult to see
that the degeneracy $D(H)$ of any graph~$H$ satisfies $m_0(H) \le D(H) \le
2m_0(H)$. The bound in the first part of the following result is better
than that in the first part of Theorem~\ref{thm:FNPuniversal} if
$D(H)<2m_0(H)-\frac12$. The bound in the second part is better than that in
Theorem~\ref{thm:CFNSuniversal} if $D(H)<(\Delta(H)-1)/2$.

\begin{theorem}[Allen, B\"ottcher, H\`an, Kohayakawa,
  Person~\cite{sparseBlowUp}]
\mbox{}\\
  For every $\eps>0$, $\Delta\ge 1$ and $D\ge 1$ there is~$C$ such that the
  random graph $G(n,p)$ a.a.s.\ is universal 
  \begin{enumerate}[label=\abc,labelindent=0pt]
    \item for all $D$-degenerate $H\in\cH(n,\Delta)$ if
      $p\ge C(\frac{\log n}{n})^{1/(2D+1)}$, and
    \item for all $D$-degenerate $H\in\cH\big((1-\eps)n,\Delta\big)$ if
      $p\ge C(\frac{\log n}{n})^{1/(2D)}$.
    \end{enumerate}
\end{theorem}

This result is a direct consequence of a sparse blow-up lemma for graphs
with bounded degeneracy (and maximum degree) established
in~\cite{sparseBlowUp}, which we shall return to in Section~\ref{sec:BUL}.


\subsection{Universality for bounded degree trees}

Recall that the result of Alon, Krivelevich and
Sudakov~\cite{AloKriSud:trees} (Theorem~\ref{thm:AKStrees}) states that
already for $p=C(\Delta,\eps)/n$ the random graph $G(n,p)$ is a.a.s.\
universal for the family $\cT\big((1-\eps)n,\Delta\big)$ of (almost
spanning) trees on $(1-\eps)n$ vertices and maximum degree at most
$\Delta$.

For spanning trees the situation is less well understood.  Hefetz,
Krivelevich, and Szab\'o~\cite{HefKriSza:trees} showed that spanning trees
with linearly long bare paths are universally a.a.s.\ contained in $G(n,p)$
for $p=(1+\eps)\log n/n$. The first universality result in $G(n,p)$ for the
entire class $\cT(n,\Delta)$ was obtained by Johannsen, Krivelevich,
Samotij~\cite{JohKriSam:universal}.  This is a consequence of the following
universality result for graphs with certain natural expansion properties. The
proof of this result relies on the embedding result of
Haxell~\cite{Haxell:FriedPipp} for large trees in graphs with suitable
expansion properties and a result of Hefetz, Krivelevich, and
Szab\'o~\cite{HefKriSza:HamCon} on Hamilton paths between any pair of
vertices in graphs with certain different expansion properties.

\begin{theorem}[Johannsen, Krivelevich, Samotij~\cite{JohKriSam:universal}]
  \label{thm:JKSuniversal}
  \mbox{}\\
  There is a constant~$c$ such that for any~$n$ and~$\Delta$ with $\log
  n\le\Delta\le c n^{1/3}$ every graph~$G$ on~$n$ vertices with
  \begin{enumerate}[label=\rom, labelindent=0pt]
    \item $|N_G(X)|\ge 7\Delta n^{2/3}|X|$ for all $X\subset V(G)$ with
      $1\le|X|<\frac{n^{1/3}}{14\Delta}$, and
    \item $e_G(X,Y)>0$ for all disjoint $X,Y\subset V(G)$ with
      $|X|=|Y|=\lceil\frac{n^{1/3}}{14\Delta}\rceil$
  \end{enumerate}
  is universal for $\cT(n,\Delta)$.
\end{theorem}

This directly implies that if $\Delta\ge\log n$ then $G(n,p)$ is a.a.s.\ 
universal for $\cT(n,\Delta)$ if $p\ge C\Delta\log n/n^{1/3}$, and hence
universality for $\cT(n,\Delta)$ with constant $\Delta$ if $p\ge C\log^2 n/n^{1/3}$.

The result of Ferber, Nenadov, and Peter~\cite{FerNenPet:universal}
discussed in the previous section improved on this when~$\Delta$ grows
slower than $n^{1/66}/(\log n)^{1/22}$. Indeed, it follows from the second part of
Theorem~\ref{thm:FNPuniversal} that $G(n,p)$ is a.a.s.\ universal for
$\cT(n,\Delta)$ if \[p\cdot n^{1/2}/(\Delta^{12}\log^3n)\to\infty\,.\]
Further, Montgomery announced in~\cite{Montgomery:DeltaT} that, using
refinements of his method for proving Theorem~\ref{thm:Montgomery},
establishing universality of $G(n,p)$ for $\cT(n,\Delta)$ with
$p=C(\Delta)\log^2 n/n$ is now within reach.

Finally, let us remark that, again, $G(n,p)$ has no chance in giving the
sparsest graph that is universal for $\cT(n,\Delta)$.  Indeed, Bhatt,
Chung, Leighton, and Rosenberg~\cite{BhaChuLeiRos} constructed $n$-vertex
graphs which are universal for $\cT(n,\Delta)$ with constant maximum degree
$C(\Delta)$. See the references in~\cite{Alon:ChaosOrder} for earlier
constructions.


\section{Resilience of random graphs}
\label{sec:robust}

In this section we study the question of how easily an adversary can
destroy copies of a graph~$H$ in
$G(n,p)$. Questions of this type date back (at least\footnote{Of course
  Tur\'an-type problems in random graphs also fall in this category and
  were studied even earlier (cf.\ Section~\ref{sec:ErdSto}).})
to~\cite{millenium} where this phenomenon was dubbed \emph{fault
  tolerance} (which also appears in~\cite{KimVu:sandwich}), but lately the term \emph{resilience} has come into vogue,
following Sudakov and Vu~\cite{SudVu:resilience}.

 
Let~$\cP$ be a monotone increasing graph property and~$\Gamma$ be a
graph. The \emph{global resilience} of~$\Gamma$ with respect to~$\cP$ is
the minimum $\eta\in\mathbb{R}$ such that deleting a suitable set of $\eta
e(\Gamma)$ edges from~$\Gamma$ results in a graph not in~$\cP$. In other
words, whenever an adversary deletes less than a $\eta$-fraction of the
edges of~$\Gamma$, the resulting graph will still be in~$\cP$.  Similarly,
in the definition of local resilience the adversary is allowed to destroy a
certain fraction of the edges incident to each vertex.  Formally, the
\emph{local resilience} of~$\Gamma$ with respect to~$\cP$ is the minimum
$\eta\in\mathbb{R}$ such that deleting a suitable set of edges, while
respecting the restriction that for every vertex $v\in V(\Gamma)$ at most
$\eta\deg_\Gamma(v)$ edges containing~$v$ are removed, results in a graph
not in~$\cP$. For $\Gamma=G(n,p)$ with $p\ge C\log n/n$ for~$C$ sufficiently
large (where we have degree concentration) this means that
for any $\eta'>\eta$ any subgraph~$G$
of~$\Gamma$ with minimum degree at least $(1-\eta')pn$ is in~$\cP$.

For the random graph~$G(n,p)$ we may then ask what is the local or global
resilience of $G(n,p)$ a.a.s.\ with respect to a property~$\cP$ for a
given~$p$? It turns out that the answer to this question usually is either
trivial, that is, basically~$0$ or~$1$, or provided by some extremal result
in dense graphs (in other words, it is as in $G(n,p)$ with $p=1$). It is
thus not surprising that the local resilience is heavily influenced by the
chromatic number of the graphs under study. To the best of my knowledge,
at present we do not know of any (subgraph) property which
does not follow the pattern just described.

\begin{question}
  Let $\pi(\cH_n)$ be the local resilience of $G(n,1)$ with respect to
  containing all graphs from $\cH_n$.  Is there any (interesting)
  family~$\cH=(\cH_n)$ of graphs such that $\pi(\cH)=\lim_{n\to\infty}\pi(\cH_n)$ exists, and the limit as~$n$ tends to infinity of
  the local resilience of $G(n,p)$ with respect to containing all graphs
 in~$\cH_n$ exists but is not in $\{0,1-\pi(\cH)\}$?
\end{question}

Let me remark that resilience and universality are orthogonal properties in
the following sense. We might ask for which probabilities $G(n,p)$ has
a.a.s.\ a certain resilience with respect to containing any fixed graph
sequence $H=(H_n)$ from a family~$\cH$, or with respect to being universal
for~$\cH$ and there is a priori no reason why the answers should turn out
the same (though we typically expect them to be). However, in contrast to
some results discussed in the previous two sections, at present the methods
available for proving resilience generally are constructive and hence allow
for universality results. On the other hand, a side effect of this is that
many of the probability bounds obtained are far from best-possible.

I will start this section with a global resilience result for small linear
sized bounded degree bipartite graphs in
Section~\ref{sec:resilience:global}, which I also use to outline one
approach often used for obtaining resilience results that relies on the
sparse regularity lemma. I then review local resilience results for cycles
in Section~\ref{sec:resilience:cycles}, for trees in
Section~\ref{sec:resilience:trees}, for triangle factors in
Section~\ref{sec:resilience:triangles}, and for graphs of low bandwidth in
Section~\ref{sec:resilience:bandwidth}.

\subsection{Global resilience}
\label{sec:resilience:global}

Obviously, any graph must have trivial global resilience with respect to
the containment of any spanning graph~$H$, since an adversary can delete
all copies of~$H$ by simply deleting all edges at some vertex. For small
linearly sized bipartite graphs~$H$, however, Alon, Capalbo, Kohayakawa,
Ruci\'nski~\cite{millenium} and Szemer\'edi, in a paper initiating research
into the area of the resilience of random graphs, proved the following
result.

\begin{theorem}[\cite{millenium}]
\label{thm:globalUniversalBip}  
  For every $\Delta\ge 2$ and $\gamma>0$ there exist $\eta>0$ and~$C$ such
  that if
  $p\ge C(\frac{\log n}{n})^{1/\Delta}$
  then $G(n,p)$ a.a.s.\ has global resilience at least $1-\gamma$ with
  respect to universality for the family $\cH(\eta n,\eta n,\Delta)$ of all
  bipartite graphs with partition classes of size $\lfloor\eta n\rfloor$
  and maximum degree at most~$\Delta$.
\end{theorem}

Note that this shows that $G(n,p)$ contains many copies of all graphs
in~$\cH(\eta n,\eta n,\Delta)$ everywhere. It is clear that such a result
cannot hold for non-bipartite~$H$ because, as any other graph, $G(n,p)$ can
be made bipartite by deleting half of its edges.  The lower bound on~$p$
though is unlikely to be optimal.

\begin{problem}
  Improve the lower bound on~$p$ in Theorem~\ref{thm:globalUniversalBip}.
\end{problem}

The proof in~\cite{millenium} of Theorem~\ref{thm:globalUniversalBip} uses
the sparse regularity lemma, which I will present and explain in more detail
in Section~\ref{sec:BUL}. The strategy is as follows. First, the sparse
regularity lemma is applied to the graph~$G$ to obtain a sparse
$\eps$-regular partition of~$V(G)$. It is then easy to show that some pair
of clusters in this partition forms a sparse $\eps$-regular pair $(V_1,V_2)$
with sufficient density. The authors of~\cite{millenium} then develop an
embedding result for bounded degree bipartite graphs with partition classes
of size $\eta'|V_1|$ and $\eta'|V_2|$ in such a pair.\footnote{This result
  is only stated for $p\ge C(\log n/n)^{1/2\Delta}$ in~\cite{millenium} though,
  for example, with the bipartite sparse blow-up lemma inferred
  in~\cite{BoeKohTar} from their techniques and from newer regularity
  inheritance results, one easily obtains from their proof the probability
  bound claimed in Theorem~\ref{thm:globalUniversalBip}.}

Most other resilience results (with the exception of the results on cycles
in the next section) mentioned in the following use proof
strategies which are variations on this basic strategy: They use the sparse
regularity lemma to obtain a regular partition, then use a result from
dense extremal graph theory on the so-called reduced graph to obtain a
suitable structure of regular pairs in this partition, and then use or
develop a suitable embedding lemma in such structures of regular pairs,
which allows one (often with substantial extra work) to embed the desired graphs.

\subsection{Local resilience for cycles}
\label{sec:resilience:cycles}

In the language of local resilience, Dirac's theorem~\cite{Dirac} states
that~$K_n$ has local resilience $1/2-o(1)$ with respect to containing a
Hamilton cycle. In this section we shall consider sparse analogues of this
result in $G(n,p)$.

Clearly, the local resilience of $G(n,p)$ with respect to containing any
graph on more than $n/2$ vertices is at most $\frac12-o(1)$, since by
deleting the edges of $G(n,p)$ in a random balanced cut we obtain a
disconnected graph with components of size at most $\frac12n$, and it can
easily be shown that each vertex loses at most $(\frac12-o(1))pn$ of its
edges. Sudakov and Vu~\cite{SudVu:resilience} then showed a
corresponding lower bound.  They proved that for every $\gamma>0$ the local
resilience of $G(n,p)$ with respect to containing a Hamilton cycle is
a.a.s.\ at least $\frac12-\gamma$ if $p>\log^4n/n$.

Smaller probabilities were first considered by Frieze,
Krivelevich~\cite{FriKri:resilience}, who proved that there are $C$ and
$\eta$ such that for $p\ge C\log n/n$ the local resilience of $G(n,p)$ for
containing a Hamilton cycle is a.a.s.\ at least $\eta$.  Ben-Shimon,
Krivelevich, Sudakov~\cite{BenKriSud:resilience1} then were able to replace
$\eta$ with $\frac16(1-\gamma)$, and then in~\cite{BenKriSud:resilience2} with
$\frac13(1-\gamma)$.  Finally Lee and Sudakov~\cite{LeeSud:resilience}
showed that also for this range of~$p$ the local resilience is
$\frac12-o(1)$.

\begin{theorem}[Lee, Sudakov~\cite{LeeSud:resilience}]
  For every $\gamma>0$ there is a constant~$C$ such that the local resilience
  of $G(n,p)$ with respect to containing a Hamilton cycle is a.a.s.\ at
  least $\frac12-\gamma$ if $p>C\log n/n$.
\end{theorem}

In~\cite{BenKriSud:resilience2} probabilities as close as possible to
the threshold for Hamiltonicity, that is, $p\ge(\log n+\log\log
n+\omega(1))/n$, are investigated, at which point the results need to be of
a different form, because vertices of degree~$2$ may exist in $G(n,p)$.
That Hamilton cycles are so well understood is connected to the fact that
with the P\'osa rotation-extension technique (see, e.g., \cite{Posa}),
which is used in the proof of all the aforementioned results, we have a
powerful tool at hand for finding Hamilton cycles.

Even smaller probabilities, where we cannot hope for Hamilton cycles any
longer, were considered by Dellamonica, Kohayakawa, Marciniszyn and
Steger~\cite{DKMS}. They show that a.a.s.\ the local resilience of $G(n,p)$ with
respect to containing a cycle of length at least $(1-\alpha)n$ is
$\frac12-o(1)$ for any $0<\alpha<\frac12$ if $p\cdot n\to\infty$.

Finally, Krivelevich, Lee and Sudakov~\cite{KriLeeSud:resilience} proved
that if $p\cdot n^{1/2}\to\infty$ then the local resilience of $G(n,p)$
with respect to being \emph{pancyclic}, that is, having cycles of all
lengths between~$3$ and~$n$, is a.a.s.\ $\frac12-o(1)$.  Here the
probability required is higher than in the results on Hamilton cycles,
which is necessary for ensuring the adversary cannot delete all triangles
(see also the remarks in Section~\ref{sec:resilience:triangles}).  An even
stronger result was proved by Lee and Samotij in~\cite{LeeSam:resilience}
who show that for the same probability a.a.s.\ every Hamiltonian subgraph
of $G(n,p)$ containing at least $\big(\frac12+o(1)\big)pn$ edges is
pancyclic.

\subsection{Local resilience for trees}
\label{sec:resilience:trees}

Koml\'os, S\'ark\"ozy, and Szemer{\'e}di~\cite{KSS:tree} showed that for every
$\gamma>0$ and every~$\Delta$ every sufficiently large $n$-vertex graph~$G$ with
minimum degree at least $(\frac12 + \gamma)n$ contains a copy of any
spanning tree~$T$ with maximum degree at most $\Delta$.
In~\cite{KSS:treeGrowing} they then extended this result to trees with
maximum degree at most $c n / \log n$. An analogue of the former result for
random graphs in the case that~$T$ is almost spanning was obtained by
Balogh, Csaba, and Samotij~\cite{BalCsaSam:resilience}. Recall that
$\cT(n,\Delta)$ is the family of all $n$-vertex trees with maximum degree
at most~$\Delta$.

\begin{theorem}[Balogh, Csaba, Samotij~\cite{BalCsaSam:resilience}]
  For all $\Delta\ge 2$ and $\gamma>0$ there is a constant~$C$ such
  that for $p\ge C/n$ the local resilience of $G(n,p)$ with respect to
  being universal for $\cT\big((1-\gamma)n,\Delta\big)$ a.a.s.\ is at least
  $\frac12-\gamma$.
\end{theorem}

The surprising aspect about this theorem is the small probability for which
it was proven to hold.  Clearly, this bound on~$p$ is sharp up to the value
of~$C$, since for smaller~$p$ the biggest component of $G(n,p)$ gets too
small to contain a tree on $\big(1-o(1)\big)n$ vertices.  Moreover, as
argued in~\cite{BalCsaSam:resilience}, at this probability we cannot ask
for, say, balanced $D(n)$-ary trees on $\big(1-o(1)\big)n$ vertices for
$D(n)\to\infty$, since we do not have enough vertices of degree $D(n)$.
Further, the factor $\frac 12$ in this result is best possible by the
discussion in the second paragraph of the previous section.  To prove their result Balogh, Csaba,
and Samotij~\cite{BalCsaSam:resilience} use an approach based on the
regularity lemma and an embedding result for trees which is a suitable
modification of the tree embedding result by Friedman and
Pippenger~\cite{FriedPipp}.

The only local resilience result for spanning trees that I am aware of
follows from Theorem~\ref{thm:resilienceBolKom} on the resilience of
$G(n,p)$ for low-bandwidth graphs, which is presented in
Section~\ref{sec:resilience:bandwidth}. It was proven by
Chung~\cite{Chung:bandwidth} that trees with constant maximum degree have
bandwidth at most $O(n/\log n)$.

\begin{theorem}[Allen, B\"ottcher, Ehrenm\"uller, Taraz~\cite{AllBoeEhrTar}]
  \label{thm:resilienceSpTrees}
  \mbox{}\\
  For all $\Delta\ge 2$ and $\gamma>0$ there is~$C$ such that for $p\ge
  C\big(\frac{\log n}{n}\big)^{1/3}$ the local resilience of $G(n,p)$ with
  respect to being universal for $\cT(n,\Delta)$ a.a.s.\ is at least
  $\frac12-\gamma$.
\end{theorem}

This probability is not believed to be optimal. Indeed, it is
conceivable that this result remains true down to the conjectured
universality threshold for $\cT(n,\Delta)$.

\begin{conjecture}
  The conclusion of Theorem~\ref{thm:resilienceSpTrees} is true for $p\ge
  C\log n/n$.
\end{conjecture}

\subsection{Local resilience for triangle factors}
\label{sec:resilience:triangles}

Corr\'adi and Hajnal~\cite{CorHaj} proved that any graph~$G$ with
$\delta(G)\ge \frac23n$ contains a triangle factor. One could then ask if
this result can be transferred to $G(n,p)$ for~$p$ sufficiently large, that
is, if the local resilience of $G(n,p)$ with respect to containing a
triangle factor a.a.s.\ is $\frac13-o(1)$.  Huang, Lee, and
Sudakov~\cite{HuaLeeSud:constp} observed that this is not the case even for
constant~$p$. Indeed, every vertex~$v$ in~$G(n,p)$ has a.a.s.\ a
neighbourhood $N(v)$ of roughly size $pn$, and every $w\in N(v)$ has
$\deg\big(w;N(v)\big)\approx p^2n$ neighbours in $N(v)$. Therefore, we can
delete all triangles containing~$v$ by removing at most roughly $p^2n<\gamma
pn$ edges at each~$w$ if~$p$ is small compared to~$\gamma$ and hence obtain a
graph without a triangle factor. With a more careful analysis it is
possible to show that we can actually choose $O(p^{-2})$ vertices and
delete all triangles containing any of these vertices by removing less that
$\gamma pn$ edges at each vertex (for the details see
\cite[Proposition~6.3]{HuaLeeSud:constp}).

So the question above should be refined to ask for an \emph{almost
  spanning} triangle factor, covering all but $O(p^{-2})$ vertices.
Balogh, Lee and Samotij~\cite{BalLeeSam:resilienceTriangles} showed that
this is indeed true if $p\ge C\big(\frac{\log n}{n}\big)^{1/2}$.
Observe that this probability is larger than the threshold $\log^{1/3}
n/n^{2/3}$ for a triangle factor as given by Theorem~\ref{thm:JKV}.
If~$p$ grows slower than $n^{-1/2}$, however, the $O(p^{-2})$ term becomes trivial.

\begin{theorem}[Balogh, Lee, Samotij~\cite{BalLeeSam:resilienceTriangles}]
  For every $\gamma>0$ there are constants~$C$ and~$D$ such that for $p\ge
  C\big(\frac{\log n}{n}\big)^{1/2}$ the local resilience of $G(n,p)$ with
  respect to the containment of an almost spanning triangle factor covering
  all but at most $Dp^{-2}$ vertices is a.a.s.\ at least $\frac13-\gamma$.
\end{theorem}

It should be remarked that a corresponding result with $Dp^{-2}$ replaced by $\eps n$
follows easily from the conjecture of Kohayakawa, {\L}uczak, and
R\"odl~\cite[Conjecture 23]{KohLucRod}, which has long been known for
triangles and was proved in full generality
in~\cite{BalMorSam,ConGowSamSch:KLR,SaxTho}. This argument will be sketched
for the purpose of illustrating the sparse regularity lemma in
Section~\ref{sec:RL}.

For proving their result Balogh, Lee,
Samotij~\cite{BalLeeSam:resilienceTriangles} develop a sparse analogue of
the blow-up lemma for the special case of triangle factors. We shall
discuss (more general) blow-up lemmas in Section~\ref{sec:BUL}.

Analogous questions concerning $H$-factors for general~$H$ were considered for
constant~$p$ in~\cite{HuaLeeSud:constp}, but the currently best bounds
follow from Theorem~\ref{thm:resilienceBolKom}, which we discuss in the
next section.

\subsection{The bandwidth theorem in random graphs}
\label{sec:resilience:bandwidth}

In~\cite{BoeSchTar:bandwidth} it was shown that for every $\Delta,r$ and
$\gamma>0$ there is $\beta>0$ such that any sufficiently large $n$-vertex
graph~$G$ with $\delta(G)\ge(\frac{r-1}{r}+\gamma)n$ contains any
$r$-colourable $H\in\cH(n,\Delta)$ with bandwidth $\bw(H)\le\beta n$. This proved a
conjecture of Bollob\'as and Koml\'os and is often referred to as the
bandwidth theorem. It is easy to argue that some restriction like the
bandwidth restriction in this result is necessary, and also that the
$\frac{r-1}{r}$ in the minimum degree is best possible; it is also known
that we cannot have $\gamma=0$ (for details see the discussions
in~\cite{BoeSchTar:bandwidth}). Further, as shown in~\cite{BoePruTarWur},
the bandwidth condition does not excessively restrict the class of
embeddable graphs. Indeed, requiring the bandwidth of a bounded degree
$n$-vertex graph to be $o(n)$ is equivalent to requiring the treewidth to
be $o(n)$ or to have no large expanding subgraphs. This implies that
bounded degree planar graphs, and more generally bounded degree
graphs defined by some (or several) forbidden minor, have bandwidth
$o(n)$.

A transference of the bandwidth theorem to $G(n,p)$ for constant~$p$ was
obtained by Huang, Lee and Sudakov~\cite{HuaLeeSud:constp}. As discussed in
the last section in such a result we cannot hope to cover all the
graphs~$H$ embedded by the bandwidth theorem. More precisely we have to ask
for at least $O(p^{-2})$ vertices of~$H$ not to be contained in a triangle.
A result for smaller~$p$ in the special case of almost spanning bipartite
graphs in $\cH\big((1-o(1)n,\Delta\big)$ with bandwidth at most $\beta n$
was obtained in~\cite{BoeKohTar} for $p\ge C(\log
n/n)^{1/\Delta}$. Recently a general sparse analogue of the bandwidth
theorem became possible with the help of the sparse blow-up lemma (see
Section~\ref{sec:BUL}). For a concise statement, let
$\cH(n,\Delta,r,\beta)$ be the class of all $r$-colourable $n$-vertex
graphs with maximum degree~$\Delta$ and bandwidth at most $\beta n$.

\begin{theorem}[Allen, B{\"o}ttcher, Ehrenm{\"u}ller, Taraz~\cite{AllBoeEhrTar}]
\label{thm:resilienceBolKom}
\mbox{}\\
   For all $\Delta,D,r$ and $\gamma>0$ there are $\beta>0$ and $C$ such that
   $(\frac1r-\gamma)$ is a.a.s.\ a lower bound on the local resilience of $G(n,p)$
   with respect to universality for all $H\in\cH(n,\Delta,r,\beta)$ such that either
   \begin{enumerate}[label=\abc,labelindent=0pt]
     \item\label{thm:resilienceBolKom:a} at least $C\max\{p^{-2},p^{-1}\log n\}$ vertices of $H$ are not in
       triangles, and $p\ge C(\log n/n)^{1/\Delta}$, or
     \item\label{thm:resilienceBolKom:b} at least $C\max\{p^{-2},p^{-1}\log n\}$ vertices of $H$ are
       in neither triangles nor~$C_4$s, and $H$ is $D$-degenerate, and 
       $p\ge C(\log n/n)^{1/(2D+1)}$.
   \end{enumerate}
\end{theorem}

Here, the term $p^{-1}\log n$ in the bound on the vertices not in
triangles is only relevant for relatively large probabilities $p>1/\log
n$. It is an artefact of our proof and we do not believe it is necessary.
Similarly, the requirement on vertices not being contained in~$C_4$
in~\ref{thm:resilienceBolKom:b} can probably be removed, but we need it for
our proof.

Observe that this theorem provides two different lower bounds on the
probability, where the second one is better if the degeneracy of~$H$ is
much smaller than its maximum degree (note though that even in this case we
require a constant bound on the maximum degree). We do not believe these
bounds to be optimal, but the bound in~\ref{thm:resilienceBolKom:a} matches
the corresponding currently known universality bound in
Theorem~\ref{thm:DKRRuniversal} and is thus well justified. Hence, the
following problem is hard.

\begin{problem}
  Improve the bounds on~$p$ in
  Theorem~\ref{thm:resilienceBolKom}.
\end{problem}

The exponent of~$n$ in~$p$ cannot be
improved beyond $1/m_2(K_{\Delta+1})=2/(\Delta+2)$.
Indeed, if~$p$ grows slower than $n^{-1/m_2(H)}$ then in $G(n,p)$
the expected number of $H$-copies containing any fixed vertex is $o(pn)$
and one can show, using a concentration inequality of Kim and
Vu~\cite{KimVu}, that in fact a.a.s.\ every vertex of $G(n,p)$ lies in at
most $\gamma pn$ copies of~$H$ (for the details see, e.g.,
\cite[Lemma~3.3]{SparseChromThr}). Hence, in this case an adversary can even
easily delete all $H$-copies without removing more than a $2\gamma$-fraction
of the edges at each vertex.

It is possible that~$2/(\Delta+2)$ is indeed the correct exponent.  A more
precise conjecture is offered in the concluding remarks
of~\cite{AllBoeEhrTar}.

A better probability bound than that in Theorem~\ref{thm:resilienceBolKom}
was very recently obtained by Noever and Steger~\cite{NoeSte:resilience}
for the special case of almost spanning squares of Hamilton cycles, which
is approximately optimal.

\begin{theorem}[Noever, Steger~\cite{NoeSte:resilience}]
  For all $\gamma>0$
  and $p\ge n^{\gamma-1/2}$ the local resilience of $G(n,p)$ with respect to
  containing the square of a cycle on at least $(1-\gamma)n$ vertices is
  a.a.s.\ at least $\frac13-\gamma$.
\end{theorem}


\section{The blow-up lemma for sparse graphs}
\label{sec:BUL}

Szemer\'edi's regularity lemma proved extremely important for much of the
progress in extremal graph theory (and other areas) over the past few
decades. Together with the blow-up lemma it also allowed for a wealth of results
on spanning substructures of dense graphs. For sparse graphs, such as
sparse random graphs or their subgraphs, the error terms appearing in the
regularity lemma though are too coarse. This inspired the development of
sparse analogues of this machinery -- which turned out to be a difficult
task.  In this section these sparse analogues are surveyed and
some very simple example applications are provided to demonstrate how they are used.
Section~\ref{sec:RL} introduces the sparse regularity lemma and explains
how it is used for obtaining resilience results.
Section~\ref{sec:Inheritance} states so-called inheritance lemmas for
sparse regular pairs, which are needed to work with the sparse blow-up
lemma. Section~\ref{sec:rgBUL} provides the sparse blow-up lemma for random
graphs, and Section~\ref{sec:BUL:apply} outlines how it is applied.

To a certain degree I assume familiarity of the reader with the dense
regularity lemma and blow-up lemma, and refer to the
surveys~\cite{Komlos:BULsurvey,KomShoSimSze:RLsurvey,KomSim:RLsurvey} for
the relevant background.

\subsection{The sparse regularity lemma}
\label{sec:RL}

In sparse versions of the regularity lemma, all edge densities are taken
relative to an ambient density~$p$. In our applications here, where we are
interested in subgraphs~$G$ of some random graph, we may always take the
edge probability of the random graph as the ambient density~$p$. In order to
state a sparse regularity lemma we need some definitions.

Let $G=(V,E)$ be a graph, and suppose~$p\in(0,1]$ and $\eps>0$ are reals.
For disjoint nonempty sets $U,W\subset V$ the \emph{$p$-density} of the
pair $(U,W)$ is defined as $d_{G,p}(U,W)=e_G(U,W)/(p|U||W|)$.  The pair
$(U,W)$ is \emph{$(\eps,d,p)$-regular} (or
\emph{$(\eps,d,p)$-lower-regular}) if there is $d'\ge d$ such that
$d_{G,p}(U',W')=d'\pm \eps$ (or if $d_{G,p}(U',W')\ge d-\eps$,
respectively) for all $U'\subset U$ and $W'\subset W$ with $|U'|\ge\eps|U|$
and $|W'|\ge\eps|W|$.  We say that $(U,W)$ is \emph{$(\eps,p)$-regular} (or
\emph{$(\eps,p)$-lower-regular}), if it is $(\eps,d,p)$-regular (or
$(\eps,d,p)$-lower-regular) for some~$d\ge d_{G,p}(U,W)-\eps$.

An \emph{$\eps$-equipartition} of~$V$ is a partition $V=V_0\dcup V_1
\dcup\dots\dcup V_r$ with $|V_0|\le\eps|V|$ and $|V_1|=\dots=|V_r|$.  An
\emph{$(\eps,p)$-regular partition} (or an \emph{$(\eps,p)$-lower-regular
  partition}) of~$G=(V,E)$ is an \emph{$\eps$-equipartition} $V_0\dcup
V_1 \dcup\dots\dcup V_r$ of $V$ such that $(V_i,V_j)$ is an
$(\eps,p)$-regular pair (or an $(\eps,p)$-lower-regular pair) in $G$ for
all but at most $\eps \binom{r}2$ pairs $ij\in\binom{[r]}2$.  The partition
classes $V_i$ with $i\in[r]$ are called the \emph{clusters} of the
partition and $V_0$ is the \emph{exceptional set}.

The sparse regularity lemma by Kohayakawa and
R\"odl~\cite{Yoshi:RL,KohRod:RL} and Scott~\cite{Scott:RL} asserts the
existence of $(\eps,p)$-regular partitions for sparse graphs $G$.
In applications of this sparse
regularity lemma one often only makes use of sufficiently dense regular
pairs in the regular partition, and the reduced graph of the partition
captures where these dense pairs are.  Formally, an $\eps$-equipartition $V_0\dcup
V_1 \dcup\dots\dcup V_r$ of a graph~$G=(V,E)$ is an
\emph{$(\eps,d,p)$-regular partition} (or \emph{$(\eps,d,p)$-lower-regular
  partition}) with \emph{reduced graph} $R$ if~$V(R)=[r]$ and the pair
$(V_i,V_j)$ is $(\eps,d,p)$-regular (or $(\eps,d,p)$-lower-regular) in $G$
whenever $ij\in E(R)$.  Observe that, given $d>0$, an
$(\eps,p)$-regular partition gives rise to an $(\eps,d,p)$-regular
partition of~$G$ with reduced graph~$R$, where~$R$ contains exactly the
edges~$ij$ such that~$(V_i,V_j)$ is $(\eps,p)$-regular
and~$d_{G,p}(V_i,V_j)\geq d-\eps$.

It then is a consequence of the sparse regularity lemma that graphs~$G$
with sufficiently large minimum degree relative to the ambient density~$p$
(and which do not have linear sized subgraphs of density much above~$p$)
allow for $(\eps,d,p)$-regular partitions with a reduced graph~$R$ of high
minimum degree. In this sense~$R$ inherits the minimum degree of~$G$. The
following lemma, which can be found e.g.\ in~\cite{AllBoeEhrTar}, makes
this precise.

\begin{lemma}[sparse regularity lemma, min.\ degree version]
  \label{lem:regularitylemma}
  \mbox{} \\
  For each $\eps >0$, $\alpha \in [0,1]$, and $r_0\geq 1$ there exists
  $r_1\geq 1$ with the following property. For any $d\in[0,1]$, any $p>0$,
  and any $n$-vertex graph $G$ with $\delta(G)\ge\alpha\cdot pn$ such
  that for any disjoint $X,Y\subset V(G)$ with
  $|X|,|Y|\ge\eps\frac{n}{r_1}$ we have
  $e(X,Y)\le(1+\frac{\eps^2}{10^3})p|X||Y|$, there is an
  $(\eps,d,p)$-regular partition of $V(G)$ with reduced
  graph~$R$ with $\delta(R) \geq
  (\alpha-d-\eps)|V(R)|$ and $r_0 \leq |V(R)| \leq r_1$.
\end{lemma}

The crucial point is that the reduced graph~$R$ in this lemma is a \emph{dense} graph,
which means that we can apply extremal graph theory results for dense
graphs to~$R$. It should be noted that analogous lemmas can easily be
formulated where other properties are inherited by the reduced graph, such
as the (relative) density of~$G$.

The regularity lemma then becomes useful in conjunction with suitable
embedding lemmas. These come in different flavours. Embedding constant
sized graphs~$H$ in systems of regular pairs in $G(n,p)$ is allowed by the
so-called \emph{counting lemma}, which even allows to give good estimates
on the number of $H$-copies. In a major breakthrough such counting lemmas
were recently established for the correct~$p$ (that is, the threshold was
established) in~\cite{BalMorSam,ConGowSamSch:KLR,SaxTho}, verifying a
conjecture of Kohayakawa, {\L}uczak, and R\"odl~\cite[Conjecture
23]{KohLucRod}. An embedding lemma for~$H$ of small linear size, on the
other hand, was provided in~\cite{ChvRand} for $p\ge C(\log
n/n)^{-1/\Delta}$. This range of~$p$ is not believed to be best
possible, but again matches the natural barrier. Finally, the blow-up lemma, which is stated in
Section~\ref{sec:rgBUL}, handles spanning graphs (for the same edge
probability~$p$).

To illustrate how the sparse regularity lemma and the embedding lemmas
interact, let us briefly sketch how to show that for $p\ge C(\log
n/n)^{1/2}$ a.a.s.\ a subgraph~$G$ of $G(n,p)$ with
$\delta(G)\ge(\frac23+\gamma)pn$ has a triangle factor
covering at least $(1-\gamma)n$ vertices for every $\gamma>0$ and~$C$
sufficiently large. Indeed, if we apply the minimum degree version of the
sparse regularity lemma (Lemma~\ref{lem:regularitylemma}) to~$G$, with
$\eps\ll d$ sufficiently small and $r_0=3$, we obtain a reduced graph~$R$
with $\delta(R)\ge(\frac23+\frac{\gamma}{2})v(R)$, which thus contains a
(spanning) triangle factor by the theorem of Corr\'adi and
Hajnal~\cite{CorHaj}. One triangle in this triangle factor corresponds to
three $(\eps,d,p)$-regular pairs in~$G$, in which, according to the sparse
counting lemma, we find one triangle. After removing the three vertices of this
triangle, what remains of the three pairs is still $(\eps',d,p)$-regular
for~$\eps'$ almost as big as~$\eps$. Hence, we can apply the counting lemma
again to find another triangle. In fact, we can repeat this process until,
say, a $\frac12\gamma$-fraction of the original three pairs is
left. Repeating this for each triangle in the triangle factor of~$R$, we
obtain an almost spanning triangle factor in~$G$ covering all but at most
the~$\eps n$ vertices of the exceptional set~$V_0$ and a
$\frac12\gamma$-fraction of $V\setminus V_0$.

\subsection{Regularity inheritance in \texorpdfstring{$G(n,p)$}{G(n,p)}}
\label{sec:Inheritance}

In the dense setting, when embedding graphs~$H$ in systems of regular pairs one often proceeds in
rounds, and for later rounds crucially relies on the following fact (and a two-sided version thereof). Assume
$(X,Y)$ is a regular pair into which we want to embed an edge $x'y'$
of~$H$. Assume further that some neighbour~$z'$ of~$x'$ was embedded in
previous rounds in a pair~$Z$. Then the setup of the blow-up lemma will be
such that $(Z,X)$ is also a regular pair, and we will have chosen the
image~$z$ of~$z'$ carefully enough so that~$z$ is ``typical'' in the pair
$(Z,X)$ in the sense that $N(z;X)$ will be of size $d|X|\gg\eps|X|$ for a
suitable constant~$d$. It then easily follows from the definition of
$\eps$-regularity that the pair $\big(N(z;X),Y\big)$ is still a regular
pair (with reduced regularity parameter), that is $\big(N(z;X),Y\big)$
\emph{inherits} regularity from $(X,Y)$. This then makes it easy to embed
the edge $x'y'$ in $(X,Y)$ such that~$x'$ is embedded into $N(z;X)$.

Trying to use a similar approach in sparse graphs we encounter the
following problem: If $(X,Y)$ and $(Y,Z)$ are $(\eps,d,p)$-regular pairs
then a ``typical'' vertex $z\in Z$ has a neighbourhood of size about $dp|X|$
in~$X$, which is much smaller than $\eps|X|$ if~$p$ goes to~$0$. Hence, it
is not clear any more that $\big(N(z;X),Y\big)$ inherits regularity from
$(X,Y)$ -- in fact, this is false in general. Fortunately, however, if we
consider regular pairs $(X,Y)$ and $(Y,Z)$ in a subgraph~$G$ of $G(n,p)$,
then it is true for most $z\in Z$ that $\big(N_G(z;X),Y\big)$ inherits
regularity from $(X,Y)$. This phenomenon was observed
in~\cite{GerKohRodSte,KohRod:pairs,ChvRand}. Based on the techniques
developed in these papers, the following regularity inheritance lemmas are
shown in~\cite{sparseBlowUp}.

\begin{lemma}[One-sided regularity inheritance~\cite{sparseBlowUp}]
    \label{lem:oneRegIn} 
    For each $\eps',d>0$ there are $\eps_0>0$ and
    $C$ such that for all $0<\eps<\eps_0$ and $0<p<1$, a.a.s.\
    $\Gamma=G(n,p)$ has the following property. 
    Let $G\subset\Gamma$ be a graph and $X,Y$ be disjoint subsets of
    $V(\Gamma)$. 
    If $(X,Y)$ is $(\eps,d,p)$-lower-regular in~$G$ and
    \[|X|\ge C\max\big(p^{-2},p^{-1}\log n\big) \quad\text{and}\quad|Y|\ge Cp^{-1}\log n\,,\]
    then the pair $\big(N_\Gamma(z;X),Y\big)$ is not
    $(\eps',d,p)$-lower-regular in~$G$ for at most $Cp^{-1}\log n$ vertices
    $z\in V(\Gamma)$.
\end{lemma}

Observe that this lemma consider neighbourhoods in
$\Gamma=G(n,p)$, rather than directly in~$G$. More specifically,
Lemma~\ref{lem:oneRegIn} establishes lower-regularity of
$\big(N_\Gamma(z;X),Y\big)$. However, since for most vertices~$z\in Z$ the
order of magnitude of $\deg_G(z;X)$ and $\deg_\Gamma(z;X)$ differs by a
factor of at most $2d$, the pair $\big(N_G(z;X),Y\big)$ then easily
inherits regularity from $\big(N_\Gamma(z;X),Y\big)$.

Lemma~\ref{lem:oneRegIn} is complemented by the following two-sided version, which
guarantees lower-regularity of the pair
$\big(N_\Gamma(z;X),N_\Gamma(z;Y)\big)$. This plays an important role when
we want to embed triangles.

\begin{lemma}[Two-sided regularity inheritance~\cite{sparseBlowUp}]
    \label{lem:twoRegIn} 
    For each $\eps',d>0$ there are $\eps_0>0$ and
    $C$ such that for all $0<\eps<\eps_0$ and $0<p<1$, a.a.s.\
    $\Gamma=G(n,p)$ has the following property. 
    Let $G\subset\Gamma$ be a graph and $X,Y$ be disjoint subsets of
    $V(\Gamma)$. 
    If $(X,Y)$ is $(\eps,d,p)$-lower-regular in~$G$ and
    \[|X|,|Y|\ge C\max\big(p^{-2},p^{-1}\log n\big)\,,\] then the pair
    $\big(N_\Gamma(z;X),N_\Gamma(z;Y)\big)$ is not
    $(\eps',d,p)$-lower-regular in~$G$ for at most $C\max(p^{-2},p^{-1}\log
    n)$ vertices $z\in V(\Gamma)$.
\end{lemma}
 
These two lemmas are similar
to~\cite[Proposition~15]{ChvRand}, and in fact equivalent when
$p=\Theta\big((\log n/n)^{1/\Delta}\big)$, but not for larger~$p$, when the
bounds on~$|X|$ and~$|Y|$ and the number of vertices~$z$ are different,
which is sometimes useful in applications.

They are moreover proved for lower-regular pairs rather than for
regular pairs (which leads to a less strong assumption, but also to a
weaker conclusion). In fact, it would be interesting to obtain analogous
lemmas for sparse regular pairs.

\begin{problem}
  Prove analogues of Lemmas~\ref{lem:oneRegIn} and~\ref{lem:twoRegIn} for
  $(\eps,d,p)$-regular pairs in subgraphs of $G(n,p)$.
\end{problem}

Lemmas~\ref{lem:oneRegIn} and~\ref{lem:twoRegIn} state that most vertices
in~$Z$ satisfy regularity inheritance properties.
In the sparse blow-up lemma, however, we will require this property from all vertices in a
cluster. 
More precisely,
let $X$, $Y$ and $Z$
be vertex sets in $G\subset\Gamma$, where~$X$ and~$Y$ are disjoint and~$X$
and~$Z$ are disjoint, but we do allow $Y=Z$. We say that $(Z,X,Y)$ has
\emph{one-sided $(\eps,d,p)$-inheritance} if for each $z \in Z$ the pair
$\big(N_\Gamma(z,X),Y\big)$ is $(\eps,d,p)$-lower-regular.  If in
addition~$X$ and~$Z$ are disjoint, then we say that $(Z,X,Y)$ has
\emph{two-sided $(\eps,d,p)$-inheritance} if for each $z\in Z$ the pair
$\big(N_\Gamma(z,X),N_\Gamma(z,Y)\big)$ is $(\eps,d,p)$-lower-regular.
When applying the sparse blow-up lemma, our approach will be to simply
remove the few vertices from each cluster whose neighbourhoods in
certain other clusters do not inherit lower-regularity (and deal with them
separately).

\subsection{The random graphs blow-up lemma}
\label{sec:rgBUL}

The purpose of this section is to state a slightly simplified version of the blow-up lemma for
random graphs proven in~\cite{sparseBlowUp}. The setup in this
blow-up lemma is as follows. We are given two graphs~$G$ and~$H$ on the
same number of vertices, where $G$ is a subgraph of the random graph
$\Gamma=G(n,p)$. The graphs~$G$ and~$H$ are endowed with partitions
$\cV=\{V_i\}_{i\in[r]}$ and $\cX=\{X_i\}_{i\in[r]}$ of their respective
vertex sets, of which we require certain properties. Firstly, the
partitions~$\cV$ and~$\cX$ need to be \emph{size-compatible}, that is,
$|V_i|=|X_i|$ for all $i\in[r]$. Secondly, $(G,\cV)$ needs to be
\emph{$\kappa$-balanced}, that is, there exists~$m$ such that $m\leq
|V_i|\leq \kappa m$ for all $i, j\in[r]$.

Further, we will have two reduced graphs~$R$ and~$R'\subset R$ on~$r$
vertices, where~$R$ represents the regular pairs of $(G,\cV)$. In fact, we
work with lower-regularity instead of regularity, because that is what the
inheritance lemmas discussed in the last section provide.  Hence, we say
that $(G,\cV)$ is an \emph{$(\eps,d,p)$-lower-regular $R$-partition} if for
each edge $ij\in R$ the pair $(V_i,V_j)$ is
$(\eps,d,p)$-lower-regular.\footnote{Observe that this differs from an
  $(\eps,d,p)$-lower-regular partition with reduced graph~$R$ in that we do
  not require the partition to be an $\eps$-equipartition. In fact, in the
  partitions referred to in the blow-up lemma the exceptional set is
  omitted.}  We require that~$H$ has edges only along lower-regular pairs of
this partition.  Formally, $(H,\cX)$ is an \emph{$R$-partition} if each
part of $\cX$ is empty, and whenever there are edges of $H$ between $X_i$
and $X_j$, the pair $ij$ is an edge of $R$.

As in the dense blow-up lemma, we cannot hope to embed a spanning graph
solely in systems of regular pairs, as these may contain isolated
vertices. Therefore, we will require certain pairs to be super-regular, that
is to additionally satisfy a minimum degree condition. Where these
super-regular pairs are is captured by the second reduced graph~$R'$.
A pair $(X,Y)$ in $G\subset\Gamma$ is called \emph{$(\eps,d,p)$-super-regular} (in~$G$) if it
is $(\eps,d,p)$-lower-regular and for every $x\in X$ and $y\in Y$ we have
\begin{align*}
  \deg_G(x;Y) &>(d-\eps) \max\{p |Y|, \deg_{\Gamma}(x;Y)/2\}\,, \\
  \deg_G(y;X) &>(d-\eps) \max\{p |X|, \deg_{\Gamma}(y;X)/2\}\,.
\end{align*}
The second term in these maxima is technically necessary to treat
vertices~$x$ of exceptionally high $\Gamma$-degree into~$Y$, but can be
ignored for most purposes (see also the discussion in~\cite{sparseBlowUp}).
The partition $(G,\cV)$ is \emph{$(\eps,d,p)$-super-regular on $R'$} if for every
$ij\in E(R')$ the pair $(V_i,V_j)$ is $(\eps,d,p)$-super-regular.

But even requiring super-regularity is not enough, as for super-regular
pairs $(X,Y)$, $(Y,Z)$, $(X,Z)$ in $G(n,p)$ there may be vertices $z\in Z$
with no edge in $\big(N_G(z;X),N_G(z;Y)\big)$, which prevents us for
example from embedding a triangle factor in $(Z,X,Y)$. However, as argued
in the previous section, lower-regularity
does not get inherited on neighbourhoods for only a few vertices in $(Z,X,Y)$. Hence, omitting these we can
circumvent this problem. In the blow-up lemma we will thus require
regularity inheritance along~$R'$.  Formally, $(G,\cV)$ has \emph{one-sided
  inheritance} on $R'$ if $(V_i,V_j,V_k)$ has one-sided
$(\eps,d,p)$-inheritance for every $ij, jk\in E(R')$, where we do allow
$i=k$.  Similarly, $(G,\cV)$ has \emph{two-sided inheritance} on $R'$ if
$(V_i,V_j,V_k)$ has two-sided $(\eps,d,p)$-inheritance
for every $ij,jk,ik\in R'$.

It remains to describe which of the edges of~$H$ are required to go along
the super-regular pairs captured by~$R'$. It turns out that we only need to
restrict a small linear fraction of the vertices of each~$X_i$ to having
their neighbours and second neighbours along~$R'$. We collect these special
vertices in a so-called buffer.  Formally, a family
$\tcX=\{\tX_i\}_{i\in[r]}$ of subsets $\tX_i\subset X_i$ is an
\emph{$(\alpha,R')$-buffer} for $H$ if for each $i\in[r]$ we have
$|\tX_i|\ge\alpha |X_i|$ and for each $x\in\tX_i$ and each $xy,yz\in E(H)$
with $y\in X_j$ and $z\in X_k$ we have $ij\in R'$ and $jk\in R'$.
The buffer sets can be chosen by the user of the blow-up lemma, which
asserts that for any graphs~$H$ and~$G$ with the setup as just described we
can embed~$H$ into~$G$ (if~$p$ is sufficiently large).

\begin{lemma}[Blow-up lemma for $G_{n,p}$~\cite{sparseBlowUp}]
\label{lem:sparseBlowUp}
  For all $\Delta\ge 2$, $\Delta_{R'}\ge 1$, $\kappa\ge 1$, and $\alpha, d>0$
  there exists $\eps>0$ such that for all~$r_1$ there is a~$C$ such that for
  $p\ge C(\log n/n)^{1/\Delta}$
  the random graph $\Gamma=G(n,p)$ a.a.s.\ 
  satisfies the following.
  Let $R$ be a graph on $r\le r_1$ vertices and let $R'\subset R$ be a spanning
  subgraph with $\Delta(R')\leq \Delta_{R'}$.
  Let $H$ and $G\subset \Gamma$ be graphs with $\kappa$-balanced
  size-compatible vertex partitions
  $\cX=\{X_i\}_{i\in[r]}$ and $\cV=\{V_i\}_{i\in[r]}$, respectively, which have
  parts of size at least $m\ge n/(\kappa r_1)$.
  Let $\tcX=\{\tX_i\}_{i\in[r]}$ be a family of subsets of $V(H)$ and
  suppose that
  \begin{enumerate}[label=\rom,labelindent=0pt]
  \item $\Delta(H)\leq \Delta$, $(H,\cX)$ is an $R$-partition, 
    and $\tcX$ an $(\alpha,R')$-buffer for~$H$,
  \item $(G,\cV)$ is an $(\eps,d,p)$-lower-regular $R$-partition, which
    is $(\eps,d,p)$-super-regular on~$R'$, and 
    has one- and two-sided inheritance on~$R'$.
  \end{enumerate}
  Then there is an embedding of~$H$ into~$G$.
\end{lemma}

A number of remarks are in place. Firstly, one of the advantages in this
formulation of the blow-up lemma, compared to that of the dense blow-up
lemma, is that the required regularity constant~$\eps$ does not depend on
the number of clusters~$r$ of the reduced graph~$R$, but only on the
maximum degree of~$R'$. This makes it possible to apply this blow-up lemma
to the whole reduced graph of a regular partition given by the sparse
regularity lemma, instead of the repeated applications of the blow-up lemma
to small parts of the reduced graph together with a technique for
``glueing'' the different so-obtained subgraphs together that were the norm
when applying the dense blow-up lemma. Since for $p=1$ we recover the dense
setting, this technique can now also be used for dense graphs~$G$.

Secondly, the version of this blow-up lemma given in~\cite{sparseBlowUp} is
stronger in the following senses. One difference is that
in~\cite{sparseBlowUp} we only require two-sided inheritance on triangles
of~$R'$ in which we want to embed some triangle of~$H$ containing a vertex
in the buffer. In particular, we do not need two-sided inheritance at all
if~$H$ has no triangles. This is useful in some applications, as explained
in~\cite{sparseBlowUp}.
The other difference is that in~\cite{sparseBlowUp} so-called image
restrictions are allowed, that is, for some vertices~$x$ of~$H$ we are allowed
to specify a relatively small set of vertices in~$G$ into which~$x$ is to
be embedded. These image restrictions have somewhat more complex
requirements than in the dense case, hence we omit them here, but the basic
philosophy is that the requirements are those needed to guarantee
compatibility with super-regularity and regularity inheritance in the
remainder of the partition of~$G$ (and they are generalisations of the
image restrictions in the dense case). 

But why do we need image restrictions? In the dense case such image
restrictions were usually used for ``glueing'' different blow-up lemma
applications together, which is now no longer needed, as described
above. However, as we will describe in the next section, when we want to
apply the blow-up lemma to a partition obtained from the sparse regularity
lemma, we will need to exclude a number of vertices from each cluster to
guarantee super-regularity and regularity inheritance. In the dense case,
usually all of these vertices can be redistributed to other clusters
without destroying these properties, but in the sparse case this is not
necessarily possible. Hence, if we want to obtain a spanning embedding
result we will first need to embed certain $H$-vertices on these
exceptional vertices of~$G$ by hand, which lead to image restrictions,
before we can apply the sparse blow-up lemma to embed the remainder of~$H$
(see, e.g., \cite{AllBoeEhrTar} for more details).

Finally, again, the lower bound on~$p$ is unlikely to
be best possible. 

\begin{problem}
  Improve the exponent of~$n$ in the lower bound on~$p$ in
  Lemma~\ref{lem:sparseBlowUp}.
\end{problem}

Even in a version of this lemma for only small linear sized graphs~$H$ this
would for example directly lead to an improvement on the known bounds on
so-called size Ramsey numbers (cf.\ \cite{ChvRand}), among many others.

Let me remark that in~\cite{sparseBlowUp} additionally a version of the blow-up
lemma for $D$-degenerate graphs~$H$ with maximum degree~$\Delta$ is
given. In this version the exponent in the power of~$n$ in the bound on~$p$
depends only on~$D$ (but the constant~$C$ still depends on~$\Delta$). Often
we can choose this exponent to be $2D+1$, and in some cases even smaller,
but the details depend on the choice of a suitable buffer and are
more involved. In applications image restrictions are often needed in
addition, which complicate the statement of a corresponding blow-up lemma
even further. It is this version of the sparse blow-up lemma which is
used to prove Theorem~\ref{thm:resilienceBolKom}\ref{thm:resilienceBolKom:b}.

\subsection{Applying the blow-up lemma}
\label{sec:BUL:apply}

As an example application of the sparse blow-up lemma presented in the last
section, let us briefly sketch how it can be used to show that for every $\gamma>0$, if $C$ is sufficiently large and $p\ge C\big(\log n/n\big)^{1/4}$,
a.a.s.\ any
subgraph~$G$ of $\Gamma=G(n,p)$ with $\delta(G)\ge(\frac12+\gamma)pn$
contains a copy of the $k \times k$ square grid $H=L_k$ with
$k=(1-\gamma)\sqrt{n}$.

We start by preparing~$G$ for the sparse blow-up lemma.  To this end, we
first apply the minimum degree version of the sparse regularity lemma
(Lemma~\ref{lem:regularitylemma}) to $G$ and obtain an $(\eps,d,p)$-regular
partition $V_0\dcup V_1\dcup\dots\dcup V_r$ with reduced graph~$R$ of
minimum degree bigger than $\frac12 v(R)$. Hence, $R$ has a Hamilton cycle
$C$, which contains a perfect matching if $v(R)$ is even (otherwise first
add one cluster to the exceptional set~$V_0$). This matching is our second
reduced graph~$R'$.  We assume without loss of generality that the Hamilton
cycle is $1,2,\dots,r$, and that the matching~$R'\subset C$ is
$\{1,2\},\{3,4\},\dots,\{r-1,r\}$.

We then have to transform the regular-partition of~$G$ into a super-regular
partition with regularity inheritance. Hence, we remove from each
cluster all those vertices violating super-regularity on~$R'$, which are at
most $\eps|V_i|$ vertices per cluster $V_i$, and all those vertices violating
(one-sided) regularity inheritance on~$R'$, which by Lemmas~\ref{lem:oneRegIn}
are at most $Cp^{-1}\log n$ vertices per cluster $V_i$. 

Next we prepare~$H$. For embedding the grid~$H$ we want to use the linear
cycle structure of the Hamilton cycle~$C$ in the reduced
graph~$R$. Therefore, let us first show that we can cut~$H$ into roughly
equal pieces along a linear structure. Indeed, any diagonal of~$H$ has at
most $\sqrt{n}$ vertices, hence by choosing appropriate diagonals as cuts
(that is, we ``cut'' along the diagonal) we can partition~$H$ into
$\frac r2$ sets $Y_1,\dots,Y_{r/2}$ of size $(2n/r)\pm\sqrt{n}$. A
$C$-partition $X_1,\dots,X_r$ of~$H$ is then obtained by letting~$X_{i-1}$
and~$X_{i}$ be the two colour-classes of $H[Y_{i/2}]$ for every even
$i\in[r]$.  Observe that most edges of~$H$ then go along the matching edges
of~$R'$.
Since vertices of the buffer $\tilde X_i$ for $i\in[r]$ should
have their first and second neighbourhood along~$R'$, we simply choose
$\alpha n$ vertices in $X_i$ as $\tilde X_i$ which are on diagonals of
distance at least~$3$ to any of the cut diagonals.

It is easy to check that we have $|V_i|\ge|X_i|$ for each $i$, so we can
add isolated vertices to each part $X_i$ of $H$ to ensure
size-compatibility, and then we can apply the sparse blow-up lemma,
Lemma~\ref{lem:sparseBlowUp}, to embed~$H$ into the remainder of~$G$.

In fact, for embedding the almost spanning~$H$ into~$G$ we could have
chosen a much simpler setup: We could have added $\alpha n$ isolated
vertices to each~$X_i$, and used these for the buffer~$\tilde X_i$. Then we
could have set $R'$ to be the empty graph. We chose to describe the more
complicated setting here though, because it is this setting which is
necessary for generalising this approach to obtain a spanning $H$-copy
in~$G$.

The idea of how to generalise the above proof is roughly as follows.  We would like to ``redistribute'' the
vertices~$v$ of~$G$ we deleted to different clusters of~$G$ where they do
not violate the required properties. Because of the minimum degree
condition on~$G$ this can easily be shown to be possible for vertices~$v$
which satisfy regularity inheritance for any pair of clusters $(X_i,X_j)$
with $ij\in R$, and which further have roughly the expected $\Gamma$-degree
in each $X_i$ with $i\in[r]$. The latter condition is necessary because of
the second term in the maximum in the definition of sparse super-regular
pairs, but it can be shown that a.a.s.\ all but at most $r\cdot Cp^{-1}\log
n$ vertices satisfy this condition. The remaining $r\cdot Cp^{-1}\log
n+r^3\cdot Cp^{-1}\log n$ vertices cannot be redistributed, and need to be
dealt with ``by hand''. The sparse blow-up lemma with image restrictions can
then be used to complete the embedding.  The details are more
complicated, because the redistribution process is iterative. In particular, we need to ensure that during the redistribution no new violations of
other vertices are created. Moreover, we also have to adapt the sizes of the clusters
of~$G$ to match the actual sizes of the~$X_i$. The details are omitted
(see~\cite{AllBoeEhrTar} for more explanations).

\section*{Acknowledgment}

I would like to thank the anonymous referee, Peter Allen, and Yury Person
for helpful comments and corrections.


\bibliographystyle{amsplain_yk}
\bibliography{LargeScale}


\myaddress

\newpage
\thispagestyle{empty}
\mbox{}
\newpage

\end{document}